\documentclass[11pt]{alf-j-l}
\usepackage{verbatim}
\usepackage{amssymb}

\def\Z{{\mathbb Z}}
\def\Q{{\mathbb Q}}
\def\R{{\mathbb R}}
\def\C{{\mathbb C}}
\def\A{{\mathcal A}}
\def\L{{\mathcal L}}
\def\O{{\mathcal O}}

\def\H{{\mathcal H}}

\def\SL{{\rm {SL}}}
\def\PSL{{\rm {PSL}}}
\def\Sp{{\rm{Sp}}}

\def\SO{{\rm SO}}
\def\Sp{{\rm Sp}}

\def\e.g.{{\it e.g.}}

\def\notdivides{{\mathchoice{\mathrel{\kern-1pt\not\!\kern4.0pt\bigm|\,}}%
{\mathrel{\kern-1pt\not\!\kern4pt\bigm|\,}}%
{\mathrel{\kern-1pt\not\!\kern3.5pt|\,}}%
{\mathrel{\kern-1pt\not\!\kern2.8pt|\,}}%
}}

\newcount\hours
\newcount\minutes
\def \SetTime{\hours=\time
\global\divide\hours by 60 \minutes=\hours \multiply\minutes by 60
\advance\minutes by-\time \global\multiply\minutes by-1 } \SetTime
\def \now{\number\hours:\ifnum\minutes<10 0\fi\number\minutes}
\def \Now{\today\ $[$\now$]$}

\newtheorem{theorem}{Theorem}[section]
\newtheorem{lemma}[theorem]{Lemma}
\newtheorem{proposition}[theorem]{Proposition}
\newtheorem*{theorem*}{Theorem}

\theoremstyle{definition}

\newtheorem*{Remark*}{Remark}
\newtheorem*{Corollary*}{Corollary}

\theoremstyle{remark}

\newtheorem*{remark*}{Remark}

\numberwithin{equation}{section}

\copyrightinfo{\number\year}{Paula B Cohen}

\begin{document}

\setcounter{page}{1}
\def\currentvolume{1}
\def\currentissue{Draft of\ }
\def\paperdate{\today}
\def\ISSN{}
\pagespan{1}{37}


\title {Hyperbolic distribution problems
\\
on Siegel 3-folds
\\
and Hilbert modular varieties\\\phantom{nothing}}

\author{Paula B. Cohen}
\curraddr{Department of Mathematics, Texas A\&M University, TAMU 3368, TX 77843-3368 
\penalty10000 USA}
\email{pcohen@math.tamu.edu}


\date{\Now.}



\begin{abstract}{We generalize to Hilbert modular varieties of arbitrary dimension the work of W.
Duke \cite{Du} on the equidistribution of Heegner points and of
primitive positively oriented closed geodesics in the Poincar\'e
upper half plane, subject to certain subconvexity results. We also
prove vanishing results for limits of cuspidal Weyl sums
associated with analogous problems for the Siegel upper half space
of degree 2. In particular, these Weyl sums are associated with
families of Humbert surfaces in Siegel 3-folds and of modular
curves in these Humbert surfaces.}
\end{abstract}

\maketitle \font\headlinefont=cmti10 scaled 800
\pagestyle{myheadings}\markboth{{\headlinefont Paula B
Cohen}}{{\headlinefont Hyperbolic equidistribution problems}}

\centerline{\emph{To appear in Duke Math J.}}

\tableofcontents

\section{Introduction}\label{s:1}

In this paper, we use the Maass correspondence for special
orthogonal groups ${\mathrm{SO}}(p,q)$, with integers $p,q\ge1$,
$p+q=m$, together with ``accidental'' isomorphisms between these
groups and certain modular groups in the case $m=3,4,5$, to derive
explicit formulae expressing cuspidal Weyl sums, associated to
hyperbolic distribution problems in Siegel 3-folds and Hilbert
modular varieties, in terms of Fourier coefficients for Maass and
Hilbert-Maass forms of half-integral weight. The case $m=3$,
$(p,q)=(2,1)$ with base field $\Q$ was studied in \cite{Du}.
Convexity and sub-convexity results for these Fourier
coefficients, combined with an analogous treatment of the
eigenfunctions for the continuous part of the spectrum of the
Laplace-Beltrami operator, imply the equidistribution properties
stated in this section.

Let $Q$ be a non-degenerate integral quadratic form whose
signature over $\R$ is $(p,q)$, $p+q=m$, $pq\not=0$. Let
$\lambda\in\R$, $\lambda\not=0$, and
$$
W_\lambda=\{x\in\R^m:Q(x)=\lambda\}.
$$
The group $G=\Omega(Q)$ of $m\times m$ matrices leaving $Q$
invariant is isomorphic to $\SO(p,q)$ and $G(\R)$ acts
transitively on $W_\lambda$. The stabilizer of any $x\in
W_\lambda$ is isomorphic to $\SO(p-1,q)$ if $\lambda>0$ and to
$\SO(p,q-1)$ if $\lambda<0$. A choice of Haar measure on $G(\R)$
determines an invariant volume form on the majorant space $\H_Q$
(see \S\ref{s:2}). Let $\Delta_Q$ be the Laplace-Beltrami operator
on $\H_Q$ and $\Gamma$ an arithmetic subgroup of $G$, given by a
congruence subgroup of a unit group of $Q$. In \cite{Ma}, Maass
constructed a $\theta$-lift on the space of $\Gamma$-invariant
$L^2$-integrable functions on $\H_Q$. This theta-lift converges on
the $\Delta_Q$-eigenfunctions for the discrete spectrum and has
image a corresponding Maass cusp form of half integral weight (see
Proposition \ref{prop:lift}).

As mentioned briefly in \cite{Du}, p84, the classical Maass
correspondence in the cases $m=4$, $m=5$ leads to the study of
other modular groups not treated in that paper. These modular
groups arise from considering families of polarized abelian
varieties of complex dimension 2. Recall that the complex points
$V(\C)$ of the Siegel 3-fold can be realized as the quotient of
the Siegel upper half space ${\H}_2$ of degree 2 by the integer
points ${\rm{Sp}}(4,\Z)$ of the symplectic group in real dimension
4,
$$V(\C)\simeq {\rm{Sp}}(4,\Z)\backslash {\H}_2.$$The
underlying variety $V$ has the structure of a quasi-projective
variety defined over $\Q$. Here,
\begin{equation}\label{upper} {\H}_2=\{z\in M_2(\C): z=z^t, {\rm{Im}}(z)> 0\}
\end{equation}
and, \begin{equation}\label{symp} {\rm Sp}(4,\R)=\{g\in
{\rm{GL}}_4(\R):g^tJg=J\},
\end{equation}
where
\begin{equation}
\label{J} J=\begin{pmatrix}0&-1_2\cr 1_2&0\cr\end{pmatrix}.
\end{equation}
The group ${\rm Sp}(4,\R)$ acts on ${\H}_2$ by
$$ z'\mapsto \gamma z=(Az+B)(Cz+D)^{-1},\qquad \gamma=\begin{pmatrix}A&B\\C&D\end{pmatrix}.$$
The projective symplectic group
$${\rm{PSp}}(4,\R)={\rm{Sp}}(4,\R)/\{\pm{\rm{Id_2}}\}$$is the
full group of analytic automorphisms of the complex domain
${\H}_2$. The matrix $J$ defines a complex structure and a
symplectic form $E$ on $\R^4$.

To every point $z\in{\H}_2$ we can associate the complex torus
$$
A=\C^2/\Z^2+z.\Z^2
$$
where $L=\Z^2+z.\Z^2$, and the Riemann form $E$ determines an
$\R$-linear non-degenerate alternating form on $\C^2\times\C^2$
taking integer values on $L\times L$ which gives a principal
polarization of $A$. The complex torus $A$ has the structure of an
abelian surface. In fact, the points of the complex variety
$V(\C)$ correspond bijectively with the complex isomorphism
classes of principally polarized abelian surfaces.

For an abelian variety $A$, we let ${\rm{End}}(A)$ be the
endomorphism ring and we put
${\mathrm{End_o}}(A)={\mathrm{End}}(A)\otimes_\Z\Q$. If $A$ is
simple, then ${\mathrm{End_o}}(A)$ is a division algebra over $\Q$
with a positive involution induced by the polarization of $A$.

Albert \cite{Alb1}, \cite{Alb2}, \cite{Alb3} classified the
division algebras over $\Q$ with positive involution. For the case
of abelian surfaces ${\rm dim}(A)=2$, this gives the following
(see for example \cite{VdG} Proposition (1.2)).

\goodbreak

\begin{proposition}\label{prop:endo}
If $A$ is a simple abelian surface then ${\rm{End}}(A)$ is one of
the following:\\
(1) the ring $\Z$,\\
(2) an order in a real quadratic field $F$,\\
(3) an order in an indefinite quaternion division algebra
${\mathcal Q}$ over $\Q$,\\
(4) an order in a quartic CM field $K$ (totally imaginary
quadratic extension of a real quadratic field).
\end{proposition}

We restrict ourselves to families of simple abelian surfaces, the
non-simple case being essentially covered by \cite{Du}, as then
the abelian surface is isogenous to a product of elliptic curves.
In Proposition \ref{prop:endo}, case (1) is the generic case which
will also not interest us here. Case (2) leads to considering
families of Hilbert modular surfaces, case (3) to families of
modular curves and case (4) to families of CM points. Case (4)
will appear as a special case of Theorem \ref{equigrh}.

The case $m=5$, $(p,q)=(3,2)$, together with a natural isomorphism
between $\Sp(4,\R)$ and $\SO(3,2)$ enables us to apply the results
of \S{\ref{s:2}} to $\H_2$ acted on by the lattice $\Sp(4,\Z)$.
Let ${\mathcal R}_1$ be a fundamental domain for the action of
$\Sp(4,\Z)$ on $\H_2$ and normalize the invariant volume form
$\omega_1$ on $\H_2$ so that ${\mathcal R}_1$ has volume
$\omega_1({\mathcal R}_1)=1$. We describe in \S\ref{s:4}
sub-domains $S_d$, $d>0$ and ${\mathcal E}_d$, $d<0$ arising from
$\lambda=d$ square-free with $d\equiv 1$ mod 4. The complex
surfaces $S_d$ are related to case (2) in that they are the
Humbert surfaces and parameterize abelian surfaces whose
endomorphism ring contains an order in the real quadratic field
$F=\Q(\sqrt{d})$. On the other hand ${\mathcal E}_d$ has the
structure of a real variety of real dimension 3.

The case $m=4$, $(p,q)=(2,2)$, together with a natural isomorphism
between $\SL(2,F)\otimes_\Q\R$, where $F=\Q(\sqrt{d})$, and
$\SO(2,2)$ leads to considering $\H^2$ acted on by the lattice
$\SL(\O\oplus\O^\vee)$, where $\O$ is the ring of integers of $F$.
Let ${\mathcal R}_2$ be a fundamental domain for the actions of
$\SL(\O\oplus\O^\vee)$ on $\H^2$ and normalize the invariant
volume form $\omega_2$ on $\H^2$ so that $\omega_2({\mathcal
R}_2)=1$. We describe in \S\ref{s:4} subdomains ${\mathcal
X}_{d,n}$, $n>0$ arising from $\lambda=n$ square-free with
$n\equiv N(\alpha)N(\O^\vee)$, some $\alpha\in F$ (the case $n<0$
leads to nothing new as $\SO(2,1)\simeq\SO(1,2)$). The complex
curves ${\mathcal X}_{d,n}$ are related to case (3) in that they
parameterize abelian surfaces whose endomorphism ring contains an
order in the indefinite quaternion algebra ${\mathcal Q}_{d,n}$
over $\Q$ with parameters $(d,-n/\delta d)$, for a certain
$\delta\in F$.

In \S\ref{s:5} we derive in Proposition \ref{weylcusplim}
vanishing results for limits of the cuspidal Weyl sums on the side
of the orthogonal groups in the above situations, which then apply
to the modular side by the discussion of \S\ref{s:4}. Although we
apply our results to families of principally polarized abelian
varieties, the same arguments go through without this polarization
assumption. In this paper, we do not explore in these same
situations the analytically more involved question of how to
modify the classical Maass correspondence for the eigenfunctions
of the continuous spectrum of $\Delta_Q$. Alternatively, one may
derive directly in the case $m=4,5$ upper bounds for Weyl sums for
eigenfunctions for the continuous spectrum, that is the analogues
of the results for $m=3$ of our \S\ref{s:7}. We hope to return to
this in a later paper. In general, the cuspidal case we treat here
is arithmetically more interesting. Together, such results give
the following.

\medskip

\noindent {\bf Equidistribution in genus 2:} {\it (i) The family
$\{S_d\}_{d>0}$ of Humbert surfaces and the family $\{{\mathcal
E}_d\}_{d<0}$ of real 3-folds, where $d$ is square-free and
$d\equiv1$ mod 4, are equidistributed in
$\Sp(4,\Z)\backslash\H_2$. Namely, if $\Omega_1$ is a convex
region with smooth boundary in ${\mathcal R}_1$ we have
\begin{eqnarray}
\label{equihum} \lim_{d\rightarrow\infty}{\frac{{\rm Vol}(S_d\cap
\Omega_1)}{{\rm Vol}(S_d)}}= \omega_1(\Omega_1),\nonumber\\
\lim_{-d\rightarrow\infty}{\frac{{\rm Vol}({\mathcal E}_d\cap
\Omega_1)}{{\rm Vol}({\mathcal E}_d)}}= \omega_1(\Omega_1),
\end{eqnarray}
(ii) Let $d$ be a positive square-free integer and $\O$ the ring
of integers of $\Q(\sqrt{d})$. The family $\{{\mathcal
X}_{d,n}\}_{n>0}$ of modular curves, where $n$ is congruent mod
$d$ to the norm of an ideal in the same class as the inverse
different $\O^\vee$ of $\Q(\sqrt{d})$, is equidistributed in
$\SL(\O\oplus\O^\vee)\backslash\H^2$. Namely, if $\Omega_2$ is a
convex region with smooth boundary in ${\mathcal R}_2$ we have
\begin{equation}
\label{equimod} \lim_{n\rightarrow\infty}{\frac{{\rm
Vol}({\mathcal X}_{d,n}\cap \Omega_2)}{{\rm Vol}({\mathcal
X}_{d,n})}}= \omega_2(\Omega_2)
\end{equation}}

\medskip
\goodbreak
The case $m=3$, $(p,q)=(2,1)$, together with a natural isomorphism
between $\SL(2,F)\otimes_\Q\R$, where $F$ is a totally real field
of degree $g$ over $\Q$, and $\SL(2,\R)^g$ leads to considering
$\H^g$ acted on by the lattice $\SL(\O\oplus\A)$ where $\A$ is a
fractional ideal in $F$. The case $g=1$ was treated in \cite{Du}.
However, for the case $g>1$, we need to adapt the classical Maass
correspondence to the Hilbert modular situation (see \S\ref{s:3})
and we need new (as yet unproven) subconvexity results (see
\S\ref{s:8}), hence our assumption of the generalized Riemann
hypothesis (GRH) in Theorem \ref{equigrh}. In order to treat the
continuous spectrum, in \S\ref{s:7} we study directly the
corresponding Eisenstein Weyl sums.

Let's recall the associated families of abelian varieties. Let $A$
be a $g$-dimensional complex torus where now $g\ge1$. Let $F$ be a
totally real field with $[F:\Q]=g$. Then $A$ has real
multiplication (RM) if ${\rm{End}}(A)$ contains an order $\O$ in
$F$. Such a complex torus always has the structure of an abelian
variety (see \cite{VdG}, p207). We assume throughout that $\O$ is
the ring of integers of $F$ and we denote the inverse different by
$\O^\vee$. Let $\A$ be a fractional ideal of $F$.

The group
$$
{\rm{SL}}(\O\oplus\A)=\left\{\begin{pmatrix}\alpha&\beta\\
\gamma&\delta\end{pmatrix}\in{\rm{SL}}(2,F): \alpha,\delta\in\O,
\beta\in\A^{-1}, \gamma\in\A\right\},
$$
acts by fractional linear transformations on $\H^g$, via the $g$
Galois embeddings of $F$ into $\R$, and induces an action of
${\rm{PSL}}(\O\oplus\A)={\rm{SL}}(\O\oplus\A)/\{\pm{\rm{Id}}_2\}$.
The quotient space
$$
{\rm{PSL}}(\O\oplus\A)\backslash\H^g
$$
is called a Hilbert modular variety and corresponds bijectively to
the complex isomorphism classes of polarized $g$-dimensional
abelian varieties $A$ with RM by ${\O}$. In particular, as a
complex torus we may write
$$
A(\C)=\C^g/(\A+z.\O)
$$
for some $z\in\H^g$ with $\A+z.\O$ embedded in $\C^g$ using the
Galois embedding of $F$ into $\R^g$.  The abelian variety $A$ is
principally polarized when $\A=\O^{\vee}$. When $g=2$, we recover
an example of case (2) in Proposition \ref{prop:endo}. There is a
natural modular embedding
$$
{\rm{PSL}}(\O\oplus\O^\vee)\backslash{\H}^2\rightarrow
{\rm{Sp}}(4,\Z)\backslash {\H}_2,
$$
which is described in detail for arbitrary $g$ in \cite{VdG},
Chapter IX, so that the Hilbert modular surfaces can be viewed as
subsurfaces of the Siegel 3-fold. In this context, they are
referred to as Humbert surfaces.

We consider in \S\S\ref{s:6},\ref{s:7},\ref{s:8} the situation
arising from $\lambda=\Delta\in\O$, $\Delta\not=0$. When $\Delta$
is totally negative this leads to the set $\Lambda_\Delta$ of
Heegner points (coming from $\SO(2)^g$) and when $\Delta$ is
totally positive to families ${\mathcal G}_\Delta$ of real
$g$-dimensional varieties (coming from $\SO(1,1)^g$), which for
$g=1$ are primitive closed geodesics (see \S\ref{s:4}). The
Heegner points correspond to abelian varieties of dimension $g$
with complex multiplication by an order in $F(\sqrt\Delta)$, and
so the case $g=2$ is related to case (4) in Proposition
\ref{prop:endo}. Let ${\mathcal R}_g$ be a fundamental domain for
the action of $\SL(2,\O)$ on $\H^g$ and normalize the invariant
volume form $\mu_g$ on $\H^g$ so that $\mu_g({\mathcal R}_g)=1$.

The main application in the case $m=3$ is as follows.

\goodbreak

\begin{theorem}\label{equigrh}

\noindent Let $F$ be a totally real number field of degree $g\ge1$
over $\Q$ with ring of integers $\O$. Assume that $F$ has class
number 1. Let $\Delta\in\O$, $\Delta\not=0$ be a generator of the
relative discriminant of $F(\sqrt\Delta)/F$. Under GRH (or rather
subconvexity), the families $\{\Lambda_\Delta\}_{\Delta\ll 0}$ of
Heegner points and $\{{\mathcal G}_\Delta\}_{\Delta\gg 0}$ of real
$g$-dimensional subvarieties are equidistributed in
$\SL(2,\O)\backslash\H^g$. Namely, if $\Omega_g$ is a region with
smooth boundary in ${\mathcal R}_g$ we have
\begin{eqnarray}
\label{equipoint} \lim_{|N(\Delta)|\rightarrow\infty, \Delta\ll 0
}\frac{{\rm Card}\left(\Lambda_\Delta\cap \Omega_g\right)}{{\rm
Card}\left(\Lambda_\Delta\right)}= \mu_g\left(\Omega_g\right),\nonumber\\
\lim_{N(\Delta)\rightarrow\infty,\Delta\gg 0}\frac{\sum_{{\mathcal
C}\in{\mathcal G}_\Delta}{\rm Vol}\left({\mathcal C}\cap
\Omega_g\right)}{{\sum_{{\mathcal C}\in{\mathcal G}_\Delta}}{\rm
Vol}\left({\mathcal C}\right)}= \mu_g(\Omega_g).
\end{eqnarray}
\end{theorem}

\goodbreak

We have made a number of simplifying assumptions which are not
essential. In order to reduce the technicalities, we have assumed
that $\A=\O$, that $F$ has class number 1 and that $\Delta$
generates a fundamental relative discriminant. The technicalities
arising from arbitrary class number and arbitrary $\A$ can be
simplified by working in the adelic language. In \cite{COU},
\cite{CU1} it was indicated how the fundamental discriminant
assumption, appearing also in \cite{Du}, can be removed for the
case $g=1$ and those same ideas may be applicable here. See also
\cite{EOh}.

The present paper is organized as follows. In \S\ref{s:2} we
recall the classical Maass correspondence of \cite{Ma} and derive
formulae in Proposition \ref{prop:weylmaass} for Fourier
coefficients of Maass forms of half-integral weight in terms of
Weyl sums. In \S\ref{s:3} we prove new results that generalize, in
Proposition \ref{prop:rhodaven} and Proposition
\ref{prop:rhodavepos}, the Maass correspondence and the Fourier
coefficient formulae to the case $m=3$, $(p,q)=(2,1)$ with base
field $F$ a totally real field of degree $g\ge1$ and arbitrary
class number. This extends results of \cite{Du}, \cite{KS},
\cite{Ma} that treat the case $g=1$.

In \S\ref{s:4} we use the ``accidental'' isomorphisms to relate
the results of \S\ref{s:2} and \S\ref{s:3} to Shimura varieties.
The case $m=4,5$ leads to studying families of Humbert surfaces in
Siegel 3-folds and of modular curves in these Humbert surfaces.
The case $m=3$ leads to studying families of Heegner points and of
certain sub-domains of real dimension $g$, which for $g=1$ are
primitive closed geodesics, in Hilbert modular varieties of
complex dimension $g$.

In \S\ref{s:5} we show how the subconvexity results of \cite{Du}
(in fact convexity results would suffice here) can be used to give
vanishing of limits of cuspidal Weyl sums.

In \S\ref{s:6} we derive in Lemma \ref{lemma:cuspn} upper bounds
for cuspidal Weyl sums in the Hilbert modular case. In
\S\ref{s:7}, we prove in Proposition \ref{prop:weyln} and
Proposition \ref{prop:weylp} new results that extend classical
formulae of Hecke \cite{He} expressing Eisenstein Weyl sums, in
the Hilbert modular case, in terms of central values of certain
$L$-series. These results are of independent interest.

The results of \S\ref{s:6} and \S\ref{s:7} combined with
subconvexity results for Fourier coefficients of Hilbert-Maass
modular forms are then used to prove Theorem \ref{equigrh}. The
corresponding subconvexity results for the holomorphic case have
been shown in \cite{CoSaPi}. We would need (in the notation of
\S\ref{s:3}, where $\Delta$ is an integer of $F$ assumed
square-free or a fundamental relative discriminant in the case of
class number 1) the Fourier coefficients $\rho(\Delta,f)$ for $f$
a cusp form with $L^2$-norm 1 or an Eisenstein series, with
eigenvalue $\lambda$ and half-integral weight $k$, to have an
upper bound in the $\Delta$-aspect as good as $\rho(\Delta,
f)\ll_{k,\epsilon}
c(\lambda)|N_{F/\Q}(\Delta)|^{-1/4-\delta+\epsilon}$ for a fixed
$\delta>0$ and a positive explicit constant $c(\lambda)$. Partial
progress towards subconvexity results in the Maass case have been
made by Gergely Harcos \cite{Har}, but the complete adaptation of
the method of \cite{CoSaPi} to the Maass case remains elusive
\footnotemark\footnotetext{As this paper goes to press, we have
learnt that Akshay Venkatesh, in work in progress \cite{Ven},
claims the required subconvexity result using other methods.}.
Such results would follow however from GRH, so our Theorem
\ref{equigrh} remains conditional. Although we do not pursue this
here, from our methods we can estimate rates of convergence in the
above equidistribution statements.

For compact maximal flats of ${\mathrm
{SL}}_n(\Z)\backslash{\mathrm {SL}}_n(\R)/{\mathrm {SO}}(n)$, an
equidistribution result has been obtained in \cite{Oh}, and this
represents a different type of equidistribution result to that of
\cite{Du}, even in the $g=1$, $n=2$ case. An equidistribution
result for Heegner points in Hilbert modular varieties using other
methods has been announced by Zhang \cite{Zh}, assuming as yet
unproven subconvexity results for Hilbert--Maass Fourier
coefficients. In \cite{Mic}, subconvexity results are obtained for
Rankin-Selberg $L$-functions which prove an equidistribution
property for incomplete orbits of Heegner points over definite
Shimura curves.

We wish to thank Peter Sarnak for suggesting to us the study of
equidistribution problems in higher genus and for many useful
discussions. We acknowledge support from the Ellentuck Fund at the
Institute for Advanced Study, Princeton, where a large part of
this research was carried out in the year 2000, as well as to
Princeton University for its hospitality in the academic year
2001/2002 during the preparation of this paper.

Just after completing the write-up of this paper, we received a
preprint of Clozel and Ullmo \cite{CU3} where similar, and more
general, equidistribution results are independently obtained. The
methods and language used in their paper are quite different, even
though in certain aspects a comparison with the present paper is
likely implicit. They use, in particular, methods in ergodic
theory due to Ratner \cite{Rat}, formulae of Waldspurger
\cite{Wal}, and generalizations of Hecke's formulae on Eisenstein
series due to Wielonsky \cite{Wie}. Their treatment of results
analogous to our Theorem \ref{equigrh} also appeals to as yet
unproven subconvexity results. In an earlier paper \cite{CU2},
these authors prove equidistribution results for certain families
of Shimura subvarieties of positive dimension
\footnotemark\footnotetext{Note added in proof: there is a sequel
to this paper by Ullmo \cite{U}. There are also two new preprints
of Zhang and Zhang-Jiang-Li \cite{Zh2}, \cite{ZhJL}.}. They use
ergodic arguments, which do not give rates of convergence, in
contrast to the methods used in the present paper.

\goodbreak

\section{The classical Maass correspondence}\label{s:2}

As in \cite{Du}, \S4, we exploit a construction of Maass forms as
integrals of certain automorphic eigenfunctions for the ring of
invariant differential operators, and in particular of the
Laplace-Beltrami operator, against Siegel theta functions. We
recall in outline this construction in order to fix notations,
referring the reader to \cite{Du}, \cite{KS}, \cite{Ma} for
details.

Let $Q$ be a symmetric $m\times m$ matrix with half-integer
off-diagonal elements and integer diagonal elements. Let $(p,q)$,
with integers $p,q\ge0$ satisfying $p+q=m$, be the signature of
$Q$. The majorant space ${\H}_Q$ of $Q$ is defined as
$$
{\H}_Q=\left\{H\in M_m(\R): H={}^tH, H>0, HQ^{-1}H=Q\right\},
$$
and is of real dimension $pq$. It is the symmetric space attached
to the group $G=\Omega(Q)$ of all real $m\times m$ matrices $g$
such that
$$
Q[g]={}^tgQg=Q,
$$
where $A[B]={}^tBAB$ for any matrices $A$, $B$ for which this
product makes sense. Indeed, the group $G$ acts transitively on
${\H}_Q$ by
$$
H\mapsto H[g],\qquad H\in{\H}_Q,\, g\in G.
$$
The isotropy group in $G$ of any $H\in {\H}_Q$ is a maximal
compact subgroup $K$. Analogous statements hold for the connected
component of the identity of $G$. An invariant metric on ${\H}_Q$
is given by,
$$
ds^2={\rm{Trace}}(H^{-1}dH H^{-1}dH).
$$
Let $\Gamma$ be any group of finite index in the unit group
$$
\Gamma_Q={\rm{SL}}(m,\Z)\cap G
$$
and let $\overline{\Gamma}$ be the quotient of $\Gamma$ by $\{\pm
{\rm Id}\}\cap\Gamma$. Then $\overline{\Gamma}$ acts
discontinuously on ${\H}_Q$ and is of finite covolume if $Q$ is
not a binary zero form, which we assume from now on. Let
$\Delta_Q$ be the Laplace-Beltrami operator on ${\H}_Q$ and let
$d\nu$ be the invariant volume measure induced by $ds^2$. Let
$\varphi=\varphi(H)$ be an eigenfunction of $\Delta_Q$ on
$\overline{\Gamma}\backslash {\H}_Q$, with eigenvalue $\lambda'$
defined by
$$
\Delta_Q\varphi+\lambda'\varphi=0.
$$
For $\varphi_1$, $\varphi_2$ functions on
$\overline{\Gamma}\backslash {\H}_Q$, define their inner product
by
$$
<\varphi_1,\varphi_2>=\frac1{{\rm{Vol}}(\overline{\Gamma}\backslash
{\H}_Q)}\int_{\overline{\Gamma}\backslash
{\H}_Q}\varphi_1{\overline{\varphi_2}}d\nu.
$$
For $z=u+iv\in{\H}$, the complex upper half plane, and
$H\in{\H}_Q$, let $R=uQ+ivH$. Following Siegel \cite{Sie2}, we
define,
$$
\theta(z)=\theta(z,H)=\sum_{h\in\Z^m}\exp(2\pi iR[h]).
$$
From its definition it follows that, for each fixed $z\in{\H}$,
the function $\theta(z,\cdot)$ on ${\H}_Q$ is left
${\overline{\Gamma}}$-invariant. Let the discriminant $D$, the
level $N$ for $Q$ and the definition of a Maass form of weight
$k$, discriminant $D$ for level $N$ be the same as in \cite{Du},
\S2, \S4. We have the following result which is Theorem 4 of
\cite{Du}, except that we use $\overline{\theta(z)}$ instead of
$\theta(z)$ (which has the effect of exchanging $p$ and $q$).

\goodbreak

\begin{proposition}\label{prop:lift}
Let $\varphi$ be a non-constant eigenfunction of $\Delta_Q$ on
$\overline{\Gamma}\backslash {\H}_Q$ with eigenvalue $\lambda'$
and $<\varphi,\varphi>$ finite. Suppose,
$$
f(z)=v^{m/4}<\varphi,{\overline{\theta(z)}}>
$$
is absolutely convergent for each $z\in{\H}$. Then $f(z)$ is a
Maass cusp form of weight $k=p-(m/2)$ and discriminant $D$ for the
congruence subgroup $\Gamma_0(N)$ of ${\rm{SL}}(2,\Z)$ and it has
eigenvalue $\lambda=\frac14(\lambda'+m-\frac{m^2}4)$.
\end{proposition}

In particular $f(z)$ will satisfy
$$
(\Delta_k+\lambda)f=0,
$$
where
$$
\Delta_k=y^2(\frac{\partial^2}{\partial
x^2}+\frac{\partial^2}{\partial y^2})-iky\frac{\partial}{\partial
x}
$$
together with a transformation rule for $\Gamma_0(N)$ with
automorphy factor depending on $k$ and $D$, and a growth
condition at the cusps.

We fix $H_0\in {\H}_Q$. As $G$ acts transitively on ${\H}_Q$, we
can write $H\in{\H}_Q$ as $H=H_0[g^{-1}]$, $g\in G$, and then,
$$
\theta(z,g):=\theta(z,H_0[g^{-1}])=\sum_{h\in\Z^m}\exp(2\pi
iuQ[h]-2\pi vH_0[g^{-1}(h)]).
$$
As a function $\theta(z,\cdot)$ on $G$, it is left
${\overline{\Gamma}}$-invariant and right $K$-invariant. Let
$\varphi(g)$ be the left $\Gamma$-invariant and right
$K$-invariant function on $G$ induced by the eigenfunction
$\varphi$ as in Proposition \ref{prop:lift}. Then, for an
appropriate choice of Haar measure $dg$ on $G(\R)$, we have
$$
f(z)=v^{m/4}\sum_{h\in\Z^m}\exp(2\pi i
uQ[h])\int_{\Gamma\backslash G}\varphi(g)\exp(-2\pi H_0[\sqrt v
g^{-1}(h)])dg.
$$
On the other hand, we can write
$$\lambda=s(1-s)=\frac14+\kappa^2,\qquad s=\frac12+i\kappa,\quad{\rm{Re}}(s)\ge\frac12.$$
We know that $f(z)$, $z=u+iv$,
$u,v\in\R$, has a Fourier expansion of the form
\begin{equation}\label{fourier}
f(z)=\rho(0)v^{\frac12+i\kappa}+\rho'(0)v^{\frac12-i\kappa}+
\sum_{d\in\Z, d\not=0}\rho(d)W_{\frac
k2{\rm{sgn}}(d),i\kappa}(4\pi |d|v)\exp(2\pi idu),
\end{equation}
where $W_{\alpha,\beta}(\cdot)$ is the classical Whittaker
function (see \cite{MagOb}; in fact for $f$ as in Proposition
\ref{prop:lift} we have $\rho(0)=\rho'(0)=0$). Therefore, for
$d\not=0$,
\begin{eqnarray}\label{coeff}
M_d(v):=&\rho(d)W_{\frac k2{\rm{sgn}}(d),i\kappa}(4\pi
|d|v)\nonumber\\=&v^{m/4}\sum_{h\in\Z^m,
Q[h]=d}\int_{\Gamma\backslash G}\varphi(g)\exp(-2\pi
H_0[g^{-1}(\sqrt v h)])dg. \end{eqnarray} For every integer
$d\not=0$, it follows from general results of \cite{Sie2} that the
number of orbits, under the action of ${\overline{\Gamma}}$, of
the solutions of $Q[h]=d$, $h\in\Z^m$, is finite. Let the
cardinality of this orbit be $H(d)$: it is a generalized class
number. Let $\{h^{(1)},\ldots,h^{(H(d))}\}$ be a set of
representatives in $\Z^m$ of these orbits and let $\Gamma_j$ be
the stabilizer of $h^{(j)}$ in $\Gamma$. Then,
\begin{eqnarray}
v^{-m/4}M_d(v)=&\int_{\Gamma\backslash G}\sum_{Q[h]=d}\exp(-2\pi
H_0[g^{-1}(\sqrt v h)])\varphi(g)dg\nonumber\\
=&\sum_{j=1}^{H(d)}\int_{\Gamma_j\backslash G}\exp(-2\pi
H_0[g^{-1}(\sqrt v h^{(j)})])\varphi(g)dg.\nonumber
\end{eqnarray}
We let,
$$
I_j=I_j(v)=\int_{\Gamma_j\backslash G}\exp(-2\pi H_0[g^{-1}(\sqrt
v h^{(j)})])\varphi(g)dg.
$$
We can compare directly with the discussion of \cite{Ma},\S5. In
terms of our notations, the notations of that paper become: $S=Q$,
$u=\varphi$, $\nu=\frac{-(m-2)}{2}+2i\kappa$, $x=2u$, $y=2v$,
$\alpha=\frac{p}{2}-\frac{m}{4}+\frac12+i\kappa$,
$\beta=\frac{q}{2}-\frac{m}{4}+\frac12+i\kappa$, and $t=d$. It is
shown there that
$$
\exp(2\pi dv)\sum_{j=1}^{H(d)}I_j(v)
$$
satisfies a second order differential equation (\cite{Ma}, (86))
and by looking at the behavior as $v\rightarrow\infty$, one sees
that it is a multiple of a standard solution of that equation
related to the Whittaker function, which fits with
(\ref{fourier}), (\ref{coeff}). Indeed, we have (\cite{Ma}, (91))
\begin{equation}\label{whitt}
M_d(v)=v^{m/4}\sum_{j=1}^{H(d)}I_j(v)=A(2\pi|d|)^{-m/4}W_{\frac
k2{\rm{sgn}}(d),i\kappa}(4\pi|d|v),
\end{equation}
for some $A\not=0$ independent of $v$.

We now describe this factor $A$. From now on $c_1$, $c_2,\ldots$
will be positive constants depending only on $Q$ and the sign of
$d$; these constants can in fact be explicitly computed. The
function $\varphi=\varphi(g)$ on $G\sim{\rm{SO}}(p,q)$ is
$K$-invariant on the right and is an eigenfunction of the
appropriately normalized Casimir operator on $G$. Fix a solution
$E$ of $Q[E]={\rm{sgn}}(d)$. We can find $l_j\in G$ such that
$$
l_j^{-1}(h^{(j)})=\sqrt{|d|} E,
$$
since $G$ acts transitively on the set
$$
\left\{x\in\R^m: Q[x]=d\not=0\right\}.
$$
We then have,
$$
I_j=\int_{\Gamma'\backslash G}\exp(-2\pi
H_0[g^{-1}(\sqrt{|d|v}E)])\varphi(l_jg)dg,
$$
where
$$
\Gamma'=l_j^{-1}\Gamma_jl_j
$$
is the stabilizer of $E$ in $G(\R)$. Let $d\gamma$ be a fixed Haar
measure on $\Gamma'(\R)$. Such a choice determines a Haar measure
$da$ on $\Gamma'(\R)\backslash G(\R)$ such that
$$
dg=d\gamma da.
$$
We have,
\begin{equation}\label{Ij}
I_j=\int_{\Gamma'(\R)\backslash G(\R)}\exp(-2\pi d v
H_0[a^{-1}(E)])\int_{\Gamma'\backslash\Gamma'(\R)}\varphi(l_j\gamma
a)d\gamma da.
\end{equation}
Let $\psi(g)=\varphi(l_jg)$; then $\psi(g)$ is also an
eigenfunction of the normalized Casimir operator with the same
eigenvalue as $\varphi(g)$. Let
$$
J_j(a)=\int_{\Gamma'\backslash\Gamma'(\R)}\psi(\gamma a)d\gamma,
$$
then
$$J_j(\gamma a k)=J_j(a),\qquad\gamma\in\Gamma'(\R), k\in K,
$$
so that $J_j(a)$ is uniquely determined by its value on
$\Gamma'(\R)\backslash G(\R)/ K$. In \cite{Ma}, this integral is
rewritten in terms of the variable $w=H_0[a^{-1}(E)]$ and is shown
to be a multiple of a standard function in $w$ by using
(\ref{Ij}) and comparing with (\ref{whitt}). Alternatively, one
may use the above discussion together with a uniqueness argument
as done in \cite{KS},(3.7) and (3.23) for the case $(p,q)=(2,1)$.
This enables us to write, for $e$ the identity of $G$,
$$
J_j(a)=J_j(e)V_\lambda(a)
$$
where $V_\lambda(a)$ is determined by the condition
$V_\lambda(e)=1$. We then have,
$$
I_j=J_j(e)\int_{\Gamma'(\R)\backslash G(\R)}\exp(-2\pi d v
H_0[a^{-1}(E)])V_\lambda(a) da.
$$
As in \cite{Ma}, (103), we can compare this directly with
(\ref{whitt}) to deduce that,
$$
\rho(d)=c_1|d|^{-m/4}\left\{\sum_{j=1}^{H(d)}\int_{\Gamma'\backslash\Gamma'(\R)}\varphi(l_j\gamma)d\gamma\right\}.
$$
We have shown the following. The notation $d\gamma$ is used again,
now to denote the induced Haar measure on $\Gamma_j(\R)$.

\begin{proposition}\label{prop:weylmaass} We have
\begin{equation}\label{rhon}
\rho(d)=c_1|d|^{-m/4}\left\{\sum_{j=1}^{H(d)}\int_{\Gamma_j\backslash\Gamma_j(\R)}\varphi(\gamma)d\gamma\right\}.
\end{equation}
\end{proposition}

We can also check this against the formula given in \cite{Ma},
pp288--289. Namely,
$$
\rho(d)=(2\pi)^{-m/4}|d|^{\frac m4-1}\alpha_d(Q,\varphi),
$$
where
$$
\alpha_d(Q,\varphi)=c_2|d|^{-\frac m2
+1}\sum_{j=1}^{H(d)}\int_{\Gamma'\backslash
\Gamma'(\R)}\varphi(l_j\gamma)d\gamma
$$
can be interpreted as Siegel's mass \cite{Sie2} of the
representation of $d$ by $Q$, weighted against $\varphi$.

\goodbreak

\section{The Maass correspondence for the Hilbert modular
case}\label{s:3}

In this section, we generalize the classical Maass correspondence
in the case $(p,q)=(2,1)$ to the Hilbert modular case. Let $F$ be
a totally real number field of degree $g$ over $\Q$. As in
\S\ref{s:1}, we let $\O$ be the ring of integers of $F$ and $\A$
be a fractional ideal of $F$. We define
$\Gamma_\A={\rm{SL}}(\O\oplus\A)$ to be the group of matrices of
determinant 1 in the maximal order in $M_2(F)$ given by
$$\begin{pmatrix}
\O&\A^{-1}\\\A&\O
\end{pmatrix}.$$
For $z=(z_1,\ldots,z_g)\in\C^g$, $z_j=u_j+\sqrt{-1}v_j$,
$u_j,v_j\in\R$, and $\alpha\in F$, let
$$
\alpha\cdot z=\alpha_1z_1+\ldots+\alpha_gz_g,
$$
with $\alpha_j$, $j=1,\ldots,g$ the Galois conjugates of $\alpha$,
and let
$$
N(v)=\prod_{j=1}^gv_j.
$$
Let $S\in M_2(\Z)^g$ have all its coordinates equal to the matrix
$$
Q=\begin{pmatrix}0&0&-2\\0&1&0\\-2&0&0\end{pmatrix}
$$
which has signature $(2,1)$. The majorant space $H_Q$ is
isomorphic to the upper half plane ${\H}$. For $z\in{\H}^g$ with
coordinates $z_j=u_j+iv_j$, $v_j>0$ and $H\in \H_Q^g$ with
coordinates $H_j\in\H_Q$, let $R$ have coordinates
$R_j=u_jQ+iv_jH_j$, $j=1,\ldots,g$. Let $\L$ be the lattice
$\A^{-1}\oplus\O\oplus\A$ in $F^3$. We define the theta function
\begin{equation}\label{thetaH}
\theta(z,H):=N(v)^{3/4}\sum_{h\in\L}\exp(2\pi i({}^th\cdot R\cdot
h)),
\end{equation}
where
$$
{}^th\cdot R\cdot h=\sum_{j=1}^g{}^th_jR_jh_j,
$$
for $h_j$, $j=1,\ldots,g$ the Galois conjugates of $h\in\L$. Let
$H_0\in H_Q^g$ have each coordinate equal to the majorant of $Q$
given by
$$
\begin{pmatrix}2&0&0\\0&1&0\\0&0&2\end{pmatrix}.
$$
Then, there are matrices $C_j$ in $\Omega(Q)$ such that
$H_j=H_0[C_j]$, $j=1,\ldots,g$. As $\Omega(Q)$ is isomorphic to
${\mathrm{SO}}(2,1)$ which is in turn isomorphic to
${\mathrm{SL}}(2,\R)$, we may write $\theta(z,H)$ as a function
$\theta(z,g)$ with $g\in G=\Omega(Q)^g$, with coordinates $g_j$,
$j=1,\ldots,g$ also viewed as elements both of ${\rm{SL}}(2,\R)$
and ${\rm{SO}}(2,1)$. To $h={}^t(h_1,h_2,h_3)\in\L$ we may
associate the matrix $${\mathbf
h}=\begin{pmatrix}h_1&h_2/2\\h_2/2&h_3\end{pmatrix}.
$$
Then $Q[h]=-4\det({\mathbf h})$ and the induced action of
${\rm{SL}}(2,\R)$ on $h$ becomes \begin{equation}\label{act}
g:{\mathbf h}\mapsto g(h):= g{\mathbf h}g^t,\qquad g\in
{\rm{SL}}(2,\R).\end{equation} Letting
$$
{\mathbf s}(x_1,x_2,x_3)=\exp(-2\pi(2x_1^2+x_2^2+2x_3^2)),
$$
we have from (\ref{thetaH})
\begin{equation}\label{thetag}
\theta(z,g)=N(v)^{3/4}\sum_{h\in\L}\exp(2\pi i
((h_2^2-4h_1h_3)\cdot u)N({\mathbf s}(\sqrt v g^{-1}(h)),
\end{equation}
where
$$
N({\mathbf s}(\sqrt{v}g^{-1}(h)))=\prod_{j=1,\ldots,g}{\mathbf
s}(\sqrt{v_j}g_j^{-1}(h_j)).
$$
By \cite{Shi4}, \S7 there is a congruence subgroup $\Gamma_1$,
and a multiplier $J$ such that for $\gamma_1\in\Gamma_1$,
$$
\theta(\gamma_1 z,g)=J(\gamma_1,z)\theta(z,g),
$$
and for $\gamma\in\Gamma_\A$, $k\in K_{\infty}=K^g$, where $K$ is
the maximal compact of ${\rm{SL}}(2,\R)$, we have
$$
\theta(z,\gamma gk)=\theta(z,g).
$$
We may adapt the discussion of \cite{KS}, \S 2 to our situation.
For $j=1,\ldots,g$, let $\Delta_{1/2}^{(j)}$ be the Laplacian in
the variable $z_j=u_j+i v_j$ given by
\begin{equation}\label{laph}
\Delta_{1/2}^{(j)}=v_j^2\left(\frac{\partial^2}{\partial
u_j^2}+\frac{\partial^2}{\partial
v_j^2}\right)-i\frac{v_j}{2}\frac{\partial}{\partial u_j}
\end{equation}
We have
\begin{equation}\label{cas}
D_{g_j}^{(j)}\theta(z,g)=4\Delta_{1/2}^{(j)}\theta(z,g)+\frac34\theta(z,g),\qquad
j=1,\ldots,g .
\end{equation}
A Maass-Hilbert form $\varphi$ may be viewed as a function on $G$
which is $K_{\infty}$-invariant on the right and which is an
eigenfunction of the Casimir operators $D_{g_j}^{(j)}$,
satisfying for $r_j\in\R$, $j=1,\ldots,g$,
$$
D_{g_j}^{(j)}\varphi(g)=\left(-\frac14-(2r_j)^2\right)\varphi(g),\qquad
r_j\in\R.
$$
Let
\begin{multline}\label{ltu}
U={\rm L}^2_{\rm{cusp}}(\Gamma_\A\backslash \H^g)
=\{\varphi:\H^g\rightarrow\C: \varphi(\gamma z)=\varphi(z),\,\gamma\in\Gamma_\A,\\
\int_{\Gamma_\A\backslash
\H^g}|\varphi|^2N(v)^{-2}N(du)N(dv)<\infty,
\int_0^1\ldots\int_0^1\varphi(x,y)N(dx)=0,\;{\rm a.e.}\;y\}.
\end{multline}
We can make $U$ into a Hilbert space using the natural inner
product. This space is invariant under the action of the unique
self-adjoint extensions of the $g$ Laplacians
\begin{equation}\label{lapz}
\Delta_0^{(j)}=v_j^2\left(\frac{\partial^2}{\partial
u_j^2}+\frac{\partial^2}{\partial v_j^2}\right),\qquad
j=1,\ldots,g,
\end{equation}
which provide a basis of the algebra of invariant differential
operators on $\H^g$. The Maass-Hilbert forms are (simultaneous)
eigenfunctions for all $g$ Laplacians $\Delta_0^{(j)}$. Let
${\mathbf{\frac12}}\in\Q^g$ be the vector with all its coordinates
equal to $\frac12$ and $\lambda=(\lambda_j)_{j=1}^g\in\C^g$. We
may write $\lambda_j=s_j(1-s_j)$, ${\rm{Re}}(s_j)\ge1/2$. A
Hilbert-Maass form $f$ for $\Gamma_1$ of weight
${\mathbf{\frac12}}$ and eigenvalue $\lambda$ is a function
$f:\H^g\rightarrow\C$ satisfying
\begin{eqnarray}
f(\gamma z)=&J(\gamma,z)f(z),\qquad \gamma\in\Gamma_1\nonumber\\
\Delta_{1/2}^{(j)}f=&\lambda_jf,\qquad j=1,\ldots,g,\nonumber
\end{eqnarray}
with polynomial growth at the cusps. Such a function of
$z=u+iv\in\H^g$, has a Fourier series development in terms of the
classical Whittaker functions of the form
\begin{equation}\label{fuv}
f(u+iv)=\rho_0(v,f)+\sum_{\xi\in\O(f,\Gamma_1), \xi\not=0
}\rho(\xi,f)N\left(W_{\frac14{\rm{sgn}}(\xi),s-\frac12}(4\pi|\xi|v)\right),
\end{equation}
where $\O(f,\Gamma_1)$ is a certain ideal in $F$ and
$$
N\left(W_{\frac14{\rm{sgn}}(\xi),s-\frac12}(4\pi|\xi|v)\right)=
\prod_{j=1}^gW_{\frac14{\rm{sgn}}(\xi_j),s_j-\frac12}(4\pi|\xi_j|v_j),
$$
with $\xi_1,\ldots,\xi_g$ the Galois conjugates of $\xi\in F$. The
form $f$ is cuspidal if it vanishes at the cusps of $\Gamma_1$.
Let
\begin{multline}
V={\rm{L}}^2_{\rm{cusp}}(\Gamma_1\backslash\H^g)=\{f:\H^g\rightarrow\C: f(\gamma z)=J(\gamma,z)f(z),\,\gamma\in\Gamma_1,\\
f\;\;{\rm cuspidal}\;\;{\rm and}\;\;{\rm square}\;\;{\rm
integrable}\;\}.
\end{multline}
Then, by exactly similar arguments to those of \cite{KS},
Proposition 2.3, we may derive the analogue of Proposition
\ref{prop:lift}.

\goodbreak

\begin{proposition}\label{prop:lifthilbert}
If $\varphi\in U$, viewed as a function on $G$, is an
eigenfunction of the $D_{g_j}^{(j)}$ with eigenvalues
$-(\frac14+(2r_j)^2)$,$j=1,\ldots,n$ then
$$
f(z)=\int_{\Gamma_\A\backslash G}\varphi(g)\theta(z,g)dg
$$
is an element of $V$ and is an eigenfunction of
$\Delta_{1/2}^{(j)}$, with eigenvalues\\ $-(\frac14+r_j^2)$,
$j=1,\ldots,g$.
\end{proposition}

Let $\Delta\in\O$ be totally negative, and let
\begin{equation}
M_{\Delta}(v)=\int_0^1\ldots\int_0^1\left(\int_{\Gamma_\A\backslash
G }\varphi(g)\theta(u+iv,g)dg\right)N\left(\exp(-2\pi i(\Delta u
)\right).
\end{equation}
Then,
\begin{equation}
M_{\Delta}(v)=N(v)^{3/4}\int_{\Gamma_\A\backslash
G}\sum_{h_2^2-4h_1h_3=\Delta}N\left({\mathbf s}(\sqrt
vg^{-1}(h))\right)\varphi(g)dg.
\end{equation}
Let $h(\Delta)$ be the number of $\Gamma_\A$-orbits of vectors
$h\in\L$ such that $h_2^2-4h_1h_3=\Delta$ and let $h^{(i)}$ be a
representative of the $i$-th orbit and $\Gamma_i$ the stabilizer
of $h^{(i)}$. Then we may write,
\begin{equation}\label{fourhil}
M_{\Delta}(v)=N(v)^{3/4}\sum_{i=1}^{h(\Delta)}\int_{\Gamma_i\backslash
G}N\left({\mathbf s}(\sqrt vg^{-1}(h^{(i)}))\right)\varphi(g)dg.
\end{equation}
When $g=1$, this corresponds to the situation considered in
\cite{KS}, (3.2). We consider the two cases $\Delta<<0$, totally
negative, and $\Delta>>0$, totally positive.

\smallskip

\noindent {\bf Case (i):} Let $\Delta$ be totally negative. Let
$h^{(1)},\ldots,h^{(h(\Delta))}$ be as above. Consider the
integral
$$
I_i=\int_GN\left({\mathbf
s}(\sqrt{v}g^{-1}(h^{(i)}))\right)\varphi(g)dg.
$$
Then (\ref{fourhil}) becomes
\begin{equation}
M_{\Delta}(v)=N(v)^{3/4}\sum_{i=1}^{h(\Delta)}\frac{1}{|\Gamma_i|}I_i.
\end{equation}
The group ${\rm{SL(2,\R)}}$ acts transitively on the $g$
hyperboloids of $x\in\R^3$ with ${}^txQx=\Delta_j$. Therefore, we
can find a
${\overline{g}^{(i)}}=({\overline{g}^{(i)}_j})_{j=1}^g\in G$ such
that
$$
({\overline{g}^{(i)}_j})^{-1}(h_j^{(i)})=\frac{\sqrt{|\Delta_j|}}
{2}\begin{pmatrix}1\\0\\1\end{pmatrix}.$$ Let
$$E=\begin{pmatrix}1\\0\\1\end{pmatrix}^g\in(\R^3)^g,
$$
and
$$\psi(g)=\varphi({\overline{g}^{(i)}}g).$$
We have
$$
I_i=\int_GN\left({\mathbf
s}\left(\sqrt{\frac{v|\Delta|}4}g^{-1}E\right)\right)\psi(g).
$$
Using as in \cite{KS}, after (3.6), the Cartan decomposition of
${\rm{SL}}(2,\R)$, we may write, as an integral over
$a=(a_j)_{j=1}^g\in\R^g$ with
$\delta(a_j)=\frac{a_j^2-a_j^{-2}}{2}$,
\begin{multline}
I_i=\int_1^{\infty}\ldots\int_1^{\infty}N\left({\mathbf
s}\left(\sqrt{\frac{v|\Delta|}{4}}\begin{pmatrix}a^{-2}&0&0\\0&1&0\\0&0&a^{-2}\end{pmatrix}E\right)\right)\times\\
\times
\left(\int_{K_\infty}\int_{K_\infty}\psi(k_1gk_2)dk_1dk_2\right)N\left(\delta(a)\right)N\left(\frac{da}{a}\right).
\end{multline}
Now $\psi(g)$ is an eigenfunction of the $D_{g_j}^{(j)}$ with the
same eigenvalues $\lambda_j$ as $\varphi$. As in \cite{KS}, we may
use uniqueness arguments to show that there is a standard
spherical function $\omega_j(g_j)$ with eigenvalue $\lambda_j$
such that $\omega_j(e)=1$ and
\begin{equation}
I_i=\varphi({\overline{g}^{(i)}})N\left(Y_\lambda\left(\sqrt{\frac{v|\Delta|}{4}}\right)\right),
\end{equation}
where
\begin{equation}
Y_{\lambda_j}(t)=\int_1^\infty{\mathbf
s}\left(t\begin{pmatrix}a_j^{-2}\\0\\a_j^2\end{pmatrix}\right)
\omega_{\lambda_j}\left(\begin{pmatrix}a_j&0\\0&a_j^{-1}\end{pmatrix}\right)\delta(a_j)\frac{da_j}{a_j}.
\end{equation}
In conclusion,
\begin{equation}
M_{\Delta}(v)=N(v)^{3/4}\sum_{i=1}^{h(\Delta)}\frac{1}{|\Gamma_i|}
\varphi({\overline{g}^{(i)}})N\left(Y_\lambda\left(\sqrt{\frac{v|\Delta|}{4}}\right)\right).
\end{equation}
From \cite{KS}, we have the asymptotic formula
$$
Y_{\lambda_j}(t)\sim\frac{\exp(-8\pi t^2)}{32\pi t^2},\qquad
t\rightarrow\infty,
$$
and therefore as $v_j\rightarrow\infty$, $j=1,\ldots,g$,
\begin{equation}\label{asymy}
N\left(Y_{\lambda}\left(\sqrt{\frac{v|\Delta|}{4}}\right)\right)
\sim\exp(-2\pi\sum_{j=1}^gv_j|\Delta_j|)\left(\prod_{j=1}^g8\pi
v_j|\Delta_j|\right)^{-1}.
\end{equation}
On the other hand, we write,
\begin{equation}\label{rhod}
M_\Delta(v)=\rho(\Delta)N\left(W_{-\frac14,ir}(4\pi|\Delta|v)\right).
\end{equation}
Then $\rho(\Delta)$ is the ``$\Delta$''-th Fourier coefficient of
the function $f(z)$ of Proposition \ref{prop:lifthilbert}.

We have the asymptotic formula as $v_j\rightarrow\infty$,
$j=1,\ldots,g$,
\begin{equation}\label{asymw}
N\left(W_{-\frac14,ir}(4\pi|\Delta|v))\right)
\sim\exp(-2\pi\sum_{j=1}^gv_j|\Delta_j|)\left(\prod_{j=1}^g4\pi
v_j|\Delta_j|\right)^{-1/4}.
\end{equation}
From equations (\ref{asymy}), (\ref{rhod}) and (\ref{asymw}) we
deduce the following result.

\goodbreak

\begin{proposition}\label{prop:rhodaven}
For $\Delta\ll0$, the ``$\Delta$''-th Fourier coefficient of the
function $f(z)$ of Proposition \ref{prop:lifthilbert} is given by,
\begin{equation}\label{rhodave}
\rho(\Delta)=2^{-g}(4\pi)^{-3g/4}|N(\Delta)|^{-3/4}\sum_{i=1}^{h(\Delta)}\frac{1}{|\Gamma_i|}\,\varphi({\overline
g}^{(i)}).
\end{equation}
\end{proposition}

\noindent{\bf Case (ii):} Let $\Delta$ be totally positive. Let
$h^{(1)}$,$\ldots$,$h^{(h(\Delta))}$ be as above. Consider the
integral
$$
I_i=\int_{\Gamma_i\backslash G}N\left({\mathbf s}(\sqrt
vg^{-1}(h^{(i)}))\right)\varphi(g)dg.
$$
Then (\ref{fourhil}) becomes
\begin{equation}
M_{\Delta}(v)=N(v)^{3/4}\sum_{i=1}^{h(\Delta)}I_i.
\end{equation}
The group ${\rm{SL(2,\R)}}$ acts transitively on the $g$
hyperboloids of $x\in\R^3$ with ${}^txQx=\Delta_j$. Therefore, we
can find an ${\ell}^{(i)}=({\ell^{(i)}_j})_{j=1}^g\in G$ such that
$$
(\ell^{(i)}_j)^{-1}(h_j^{(i)})=\sqrt{\Delta_j}\begin{pmatrix}0\\1\\0\end{pmatrix}.$$
Let
$$
E'=\begin{pmatrix}0\\1\\0\end{pmatrix}^g\in\left(\R^3\right)^g.
$$
and
$$\psi(g)=\varphi(\ell^{(i)}g).$$
We have
$$
I_i=\int_{\Gamma_i'\backslash G}N\left({\mathbf
s}\left(\sqrt{v\Delta}g^{-1}E'\right)\right)\psi(g)dg,
$$
where
$$
\Gamma_i'={\ell^{(i)}}^{-1}\Gamma_i\ell^{(i)}.
$$
Suppose from now on that $\Delta$ is not a square. The group
$\Gamma'(\R)=\Gamma_i'(\R)$ is the $g$-th power of the stabilizer
of $\begin{pmatrix}0\\1\\0\end{pmatrix}$, and can be written as
$$
\Gamma'(\R)=\prod_{j=1}^g\left\{\pm\begin{pmatrix}p_j^{1/2}&0\\0&
p_j^{-1/2}
\end{pmatrix},\qquad 0<p_j<\infty\right\}.
$$
The group $\Gamma_i'$ is a discrete free abelian subgroup of
$\Gamma'(\R)$ of rank $g$ over $\Z$, see for example \cite{Efr1},
Chapter 1, Section 5. We can decompose each component of
$g=(g_j)_{j=1}^g\in G$ as
\begin{eqnarray}
{}&g_j=\begin{pmatrix}1&\xi_j\\0&1\end{pmatrix}\begin{pmatrix}p_j^{1/2}&0\\0&p_j^{-1/2}\end{pmatrix}k_j\nonumber\\
=&\begin{pmatrix}p_j^{1/2}&0\\0&p_j^{-1/2}\end{pmatrix}\begin{pmatrix}1&\xi_j/p_j\\0&1\end{pmatrix}k_j,
\qquad k_j\in K,\;0<p_j<\infty,\;-\infty<\xi_j<\infty.\nonumber
\end{eqnarray}
We have
$$
g_j^{-1}\begin{pmatrix}0\\1\\0\end{pmatrix}=\begin{pmatrix}-\xi_j/p_j\\1\\0\end{pmatrix}.
$$
Let $t_j=\xi_j/p_j$, then $N(dt):=\prod_{j=1}^g dt_j$ is a Haar
measure on $\Gamma'(\R)\backslash G(\R)$ and
$N(\frac{dp}p):=\prod_{j=1}^g\frac{dp_j}{p_j}$ is a Haar measure
on $\Gamma'(\R)$. We may assume that $dg=N(dt)N(\frac{dp}{p})$.

With these notations, we may write
\begin{multline}
I_i=\int_{-\infty}^{\infty}
\ldots\int_{-\infty}^{\infty}\exp\left(-2\pi\sum_{j=1}^gv_j\Delta_j(2t_j^2+1)\right)\times\\
\times\int_{\Gamma_i'\backslash\Gamma'(\R)}\psi\left(\begin{pmatrix}p^{1/2}
&0\\0&p^{-1/2}\end{pmatrix}\begin{pmatrix}1
&t\\0&1\end{pmatrix}\right)N\left(\frac{dp}{p}\right)
N\left(dt\right).
\end{multline}
Let
$$
J_i(e)=\int_{\Gamma_i'\backslash\Gamma'(\R)}\psi\left(\begin{pmatrix}p^{1/2}
&0\\0&p^{-1/2}\end{pmatrix}\right)N\left(\frac{dp}{p}\right)
=\int_{\Gamma_i\backslash\Gamma_i(\R)}\varphi\left(\gamma\right)d\gamma,
$$
where $d\gamma$ is the invariant measure on $\Gamma_i(\R)$ induced
by $N\left(\frac{dp}{p}\right)$.

Arguing as in \cite{KS}, Case (ii) (where $g=1$), we can again use
the Casimir operators to see that there is, for each
$j=1,\ldots,g$ a unique even function $V_{\lambda_j}(t_j)$ of
$t_j\in\R$ determined by the condition $V_{\lambda_j}(0)=1$ and
such that
\begin{multline}\label{Ii}
I_i=J_i(e)\int_{-\infty}^{\infty}
\ldots\int_{-\infty}^{\infty}\exp\left(-2\pi\sum_{j=1}^gv_j\Delta_j(2t_j^2+1)\right)\times\\
\times V_{\lambda_1}(t_1)\ldots V_{\lambda_g}(t_g)dt_1\ldots dt_g.
\end{multline}
From the asymptotics in \cite{KS} we have, as
$v_j\rightarrow\infty$,$j=1,\ldots,g$,
$$
I_i\sim
J_i(e)2^{-g}\prod_{i=1}^g\exp\left(-2\pi\sum_{j=1}^gv_j\Delta_j\right)\left(v_j\Delta_j\right)^{-1/2}.
$$
Hence,
\begin{equation}\label{MDvp}
M_{\Delta}(v)\sim
N(v)^{1/4}N(\Delta)^{-1/2}2^{-g}\prod_{i=1}^g\exp\left(-2\pi\sum_{j=1}^gv_j\Delta_j\right)\sum_{i=1}^{h(\Delta)}J_i(e).
\end{equation}
We also have the asymptotic formula for $j=1,\ldots,g$,
$$
W_{\frac14, ir_j}(4\pi\Delta_jv_j)\sim\exp(-2\pi
v_j\Delta_j)\left(4\pi\Delta_jv_j\right)^{1/4},\qquad
v_j\rightarrow\infty.
$$
On the other hand, we write,
\begin{equation}\label{rhodp}
M_\Delta(v)=\rho(\Delta)N\left(W_{\frac14, ir}(4\pi\Delta
v)\right).
\end{equation}
Then $\rho(\Delta)$ is the ``$\Delta$''-th Fourier coefficient of
the function $f(z)$ of Proposition \ref{prop:lifthilbert}.

We have the asymptotic formula for $v_j\rightarrow\infty$ and
$j=1,\ldots,g$,
\begin{equation}\label{MDrho}
M_\Delta(v)\sim(4\pi)^{g/4}\rho(\Delta)\exp\left(-2\pi\sum_{j=1}^gv_j\Delta_j\right)N(\Delta)^{1/4}N(v)^{1/4}.
\end{equation}
From equations (\ref{MDvp}), (\ref{rhodp}), (\ref{MDrho}) we
deduce the following result.

\goodbreak

\begin{proposition}\label{prop:rhodavepos}
For $\Delta\gg0$, the ``$\Delta$''-th Fourier coefficient of the
function $f(z)$ of Proposition \ref{prop:lifthilbert} is given by,
\begin{equation}\label{rhodavep}
\rho(\Delta)=
2^{-g}(4\pi)^{-g/4}N(\Delta)^{-3/4}\sum_{i=1}^{h(\Delta)}\int_{\Gamma_i\backslash\Gamma_i(\R)}\varphi\left(\gamma\right)d\gamma.
\end{equation}
\end{proposition}

\goodbreak

\section{Families of symmetric domains}\label{s:4}

We describe the families of symmetric domains to which we will
apply the Maass correspondence of \S\ref{s:2} and \S\ref{s:3}.
These will correspond in particular to subvarieties of the Siegel
modular variety of genus 2 and to certain Heegner points in
arbitrary genus, coming from Hilbert modular varieties.

We exploit a natural isomorphism between $\SO(3,2)$ and
$\Sp(4,\R)$, following \cite{Sie1}, X. Let $Q$ be the quadratic
form on $\R^5$ of signature $(3,2)$ given by
\begin{equation}
\label{Q} Q(x)=Q(x_1,\ldots,x_5)=x_2^2-4x_3x_1-4x_4x_5,\qquad
x={}^t(x_1,\ldots,x_5)\in\R^5.
\end{equation}
Then $\H_2$ is isomorphic to the space of vectors
$Z={}^t(z_1,\ldots,z_5)\in\C^5$ with $z_5=1$ and
$$
Z^tQZ=0,\qquad {\overline{Z}}^tQZ<0,\qquad{\rm{Im}}(z_1)> 0.
$$
We recover the description of \S\ref{s:1} by setting
$$z=\begin{pmatrix}z_1&z_2\cr z_2&z_3\end{pmatrix},$$
the lower Siegel half space of degree 2 corresponding to the
condition ${\rm{Im}}(z_1)< 0$.

We may also introduce $Q$ as the quinary quadratic form on the
space $V$ of alternating matrices of the form $M(x)$, $x\in\R^5$,
where
$$
M(x)=\begin{pmatrix}0&-2x_4&x_2&-2x_1\cr 2x_4&0&2x_3&-x_2\cr
-x_2&-2x_3&0&2x_5\cr 2x_1&x_2&-2x_5&0\cr\end{pmatrix}.
$$
With $J$ the standard symplectic $4\times4$ matrix as in
\S\ref{s:1}, we have
$$
{}^tM(x)JM(x)=Q(x)\cdot J,\qquad x\in\R^5,
$$
and this defines $(V,Q)$. The isomorphism between $\SO(3,2)$ and
$\Sp(4,\R)$ can then be seen via the action of $g\in\Sp(4,\R)$ on
$M(x)\in V$ preserving $Q$ and given by
$$
g:M(x)\rightarrow gM(x){}^tg.
$$
Fix $\lambda\in\R$, $\lambda\not=0$. The group $\SO(Q)$ acts
transitively on the solutions $x\in\R^5$ of $Q(x)=\lambda$ and the
isotropy group of any such $x$ is isomorphic to $\SO(2,2)$ if
$\lambda>0$ and to $\SO(1,3)$ if $\lambda<0$. For $Q(x)>0$, let
$$
{\mathcal R}_x=\{z\in\H_2:\begin{pmatrix}z&1_2\end{pmatrix}
M(x){}^t\begin{pmatrix}z&1_2\end{pmatrix}=0\}.
$$
For $Q(x)<0$, let
$$
{\mathcal R}^-_x=\{z\in\H_2:\begin{pmatrix}{\overline
z}&1_2\end{pmatrix} M(x){}^t\begin{pmatrix}z&1_2\end{pmatrix}=0\}
$$
By checking at $z=\sqrt{-1}I_2\in\H_2$ and using transitivity one
sees that the domains ${\mathcal R}_x$, $Q(x)>0$ are real
isomorphic to the symmetric space for $\SO_o(2,2)$ and complex
isomorphic to $\H^2$. On the other hand, the domains ${\mathcal
R}^-_x$, $Q(x)<0$, are real isomorphic to the symmetric space for
$\SO_o(3,1)$. For $d\in\Z$ let
$$
W_d=\{h\in\Z^5: Q(h)=d\}.
$$
Assume from now on that $d$ is square-free. When $d>0$, let $S_d$
be the complex surface in $\Sp(4,\Z)\backslash \H_2$ given by the
union of the images of the ${\mathcal R}_h$, $h\in W_d$ (with $h$
primitive as $d$ is square-free). The surface $S_d$ will be
non-trivial if and only if $d\equiv1$ mod 4. Then $S_d$ is called
the Humbert surface of invariant $d$. For a general reference on
Humbert surfaces see \cite{VdG}, Chapter IX. The components of the
surface $S_d$ are images of Hilbert modular surfaces, induced by
the identification of the $\O$-module $\O\oplus\O^{\vee}$ (with
the standard alternating form derived from the trace of
$F=\Q(\sqrt{d})$ over $\Q$) with the $\Z$-lattice $\Z^4$ (with the
standard symplectic form). This amounts to viewing Hilbert modular
surfaces as sub-varieties of Siegel 3-folds. In particular, by
\cite{VdG}, Chapter IX, Proposition (2.3), the abelian surface
$$
A(\C)=\C^4/\Z^4+z\cdot\Z^4
$$
has endomorphism ring ${\rm{End}}(A)$ containing $\O$ if and only
if $z$ mod $\Sp(4,\Z)$ is in $S_d$.

When $d<0$, let ${\mathcal E}_d$ be the real 3-dimensional variety
in $\Sp(4,\Z)\backslash \H_2$ given by the union of the images of
the ${\mathcal R}^-_h$, $h\in W_d$ (with $h$ primitive as $d$ is
square-free). We have ${\mathcal E}_d$ non-trivial if and only if
$-d\equiv3$ mod 4.

We now turn to studying the Hilbert modular surfaces $X_d$ with
$$
X_d(\C)\simeq\PSL(\O\oplus\O^{\vee})\backslash\H^2
$$
where $\O$ is the ring of integers of $F=\Q(\sqrt{d})$, $d>0$
square-free. There is an isomorphism between $\SO(2,2)$ and
$\SL(2,F)\otimes_\Q\R$. Fix an integral ideal $\A$ in the same
genus as $\O^{\vee}$ and let $\delta=N(\A)$ be the norm of $\A$.
Let $\sigma$ be the non-trivial Galois automorphism of $F$
determined by $\sigma:\sqrt{d}\mapsto-\sqrt{d}$. As in \cite{Kud1}
and in \cite{VdG}, Chapter V (but with some minor differences in
conventions) we let
$$
Y_d=\{M\in M_2(F):
M=\begin{pmatrix}a\sqrt{d}&\alpha\\-\alpha^\sigma&
b\sqrt{d}\end{pmatrix},\,\alpha\in F\,, a,b\in\Q\}.
$$
As a $\Q$-vector space $Y_d$ is isomorphic to $\Q^4$. Define
$Q_d:Y_d\rightarrow\Q$ to be the quadratic form given by
$$
Q_d[M]=\det M=abd+\alpha\alpha^\sigma.
$$
Then $Q_d$ has signature $(2,2)$ and we may embed $\SL(2,F)$ into
$\SO(Q)$ by the action
$$
g:M\mapsto g^{\sigma}Mg^{-1},\qquad g\in\SL(2,F).
$$
This induces a representation
$$
\SL(2,F)\otimes_\Q\R\simeq\SL_2(\R)\times\SL_2(\R)\rightarrow\SO(Q)\simeq\SO(2,2)
$$
and a corresponding isomorphism between $\H^2$ and the majorant
space $\H_d=\H_{Q_d}$ of $Q_d$. In $Y_d$ we can define the
lattice of ``integral elements'' given by
$$
Y_d(\Z)=
\{M=\begin{pmatrix}a\sqrt{d}&\alpha\\-\alpha^{\sigma}&b\sqrt{d}/\delta\end{pmatrix}:\,\alpha\in\A^{-1},\,a,b\in\Z\}.
$$
For $\lambda\in\R$, $\lambda\not=0$ the group $\SO(Q_d)$ acts
transitively on the $M\in Y_d$ with $Q_d[M]=\lambda$ and the
isotropy group of any such $M$ is isomorphic to $\SO(2,1)$.
Therefore, one may assume that $\lambda> 0$.

For $M\in Y_d$ with $Q_d[M]>0$ let
$$
\H_M=\{(z_1,z_2)\in\H^2:\,\begin{pmatrix}z_2,1\end{pmatrix}M\begin{pmatrix}z_1\\1\end{pmatrix}=0\}.
$$
By checking at $z=(\sqrt{-1},\sqrt{-1})$ and using transitivity,
we see that $\H_M$ is isomorphic to the symmetric space for
$\SO(2,1)\sim\SL(2,\R)$. Namely, it is the graph of a fractional
linear fractional transformation and is therefore a copy of $\H$
embedded into $\H^2$.

Let $n$ be a positive square-free integer. For $M\in Y_d(\Z)$
with $Q_d[M]=n$ let
$$
\Gamma_M=\{g\in\SL(\O\oplus\O^{\vee}):\,{}^tg^\sigma Mg=M\}.
$$
Let ${\mathcal X}_M$ be the image of the curve
$\Gamma_M\backslash\H_M$ in
$\SL(\O\oplus\O^{\vee})\backslash\H^2$. Finally ${\mathcal
X}_{d,n}$ is defined as the curve given by the union of all the
${\mathcal X}_M$, $M\in Y_d(\Z)$, $Q_d[M]=n$. It is called a
modular curve and is non-trivial if and only if for some
$\alpha\in F$,
$$
n\equiv N(\alpha)N(\A)\quad {\rm{mod}}\, d.
$$
Moreover, from \cite{VdG}, p102, all irreducible components of
${\mathcal X}_{d,n}$ have the same volume. The curve ${\mathcal
X}_{d,n}$ corresponds to abelian surfaces whose endomorphism ring
contains an order in a quaternion algebra. Namely, let ${\mathcal
Q}_{d,n}$ be the quaternion algebra over $\Q$ with parameters
$(d, -n/\delta d)$: it has basis elements $1$, $i$, $j$, $k$ where
$$
i^2=d,\quad j^2=-n/\delta d,\quad k=ij=-ji.
$$
For $M\in Y_d(\Z)$ with $Q_d[M]=n$ the following algebra is
isomorphic to ${\mathcal Q}_{d,n}$ (see \cite{VdG}, Chapter V,
Proposition (1.5)),
$$
{\mathcal Q}_M=\{g\in M_2(F):\,{}^tg^\sigma Mg=\det(g)M\}
$$
and contains the order of discriminant $n^2$ given by
$$
\O_M={\mathcal Q}_M\cap\begin{pmatrix}\O&\A^{-1}\\
\A&\O\end{pmatrix}.
$$
For an abelian surface $A$, we have ${\rm{End}}(A)$ contains
$\O_M$ if and only if
$$
A(\C)\simeq\C^2/\O^\vee+z.\O
$$
with $z=(z_1,z_2)\in\H_M$.

We now turn to studying the case treated in \S\ref{s:3} where
$m=3$, $(p,q)=(2,1)$ but we work over a totally real field $F$,
with $[F:\Q]=g$ and so
$$
\SL(2,F)\otimes_\Q\R\simeq\SL(2,\R)^g\simeq\SO(2,1)^g.
$$
Let $\sigma_1,\ldots,\sigma_g$ be the Galois embeddings of $F$
into $\R$ and for $\alpha\in F$, let $\alpha^{(j)}$,
$j=1,\ldots,g$ denote its Galois conjugates. Let $\A$ be a
fractional ideal of $F$ and let $\L=\A^{-1}\oplus\O\oplus\A$ in
$F^3$. Let $Q:\R^3\rightarrow\R$ be the quadratic form
$$
Q(x)=x_2^2-4x_1x_3,\qquad x={}^t(x_1,x_2,x_3)\in\R^3.
$$
Let $\Delta\in\O$, $\Delta\not=0$, and consider the set
$$
W_\Delta=\{h\in\L:\,Q(h)=\Delta\}.
$$
As in \S\ref{s:3}, let $h(\Delta)$ be the number of
$\Gamma_\A=\SL(\O\oplus\A)$-orbits of vectors $h\in W_\Delta$. Let
$h=(\alpha,\beta,\gamma)$ and consider the quadratic equations,
\begin{equation}\label{quad}
\alpha^{(j)}z^2+\beta^{(j)}z+\gamma^{(j)}=0,\qquad j=1,\ldots,g.
\end{equation}
Then, if $\alpha\not=0$, we associate to these equations the
points
$$
z_j^{\pm}=\frac{-\beta^{(j)}\pm\sqrt{\Delta^{(j)}}}{2\alpha^{(j)}},\qquad
j=1,\ldots,g,
$$
where on the right hand side of (\ref{quad}) we {\emph{choose}}
$z_j^{+}\in\H$ if $\Delta\ll0$ and we {\emph{choose}}
$z_j^+>z_j^-$ if $\Delta\gg0$ (these choices will depend on the
sign of $\alpha^{(j)}$). If $\alpha=0$, let
$z_j^+=\sqrt{-1}\infty$ and $z_j^-=\gamma_j/\beta_j$. Let
$z_h=(z_j^+)_{j=1}^g\in\C^g$ and let $\Gamma_h$ be the stabilizer
in $\Gamma_\A$ of $z_h^+$.

If $\Delta\ll0$ (totally negative) let,
$$
\Lambda_{\Delta}=\{z_h=(z_j^+)_{j=1}^g\in\H^g\;{\rm{mod}}\;\Gamma_\A,\,
h\in W_\Delta\}.
$$
We refer to this as the set of Heegner points associated to
$\Delta$. It has cardinality $h(\Delta)$.

If $\Delta\gg0$ (totally positive), then $\Gamma_h$ is a discrete
subgroup of rank $g$ of $\Gamma_h(\R)\simeq(\R_{>0})^g$, embedded
into $G(\R)=\PSL(2,\R)^g$ using Galois embeddings. The image of
$\Gamma_h\backslash\Gamma_h(\R)$ in $\Gamma\backslash G(\R)$
determines $C_h$, the real $g$ dimensional variety in
$\Gamma_\A\backslash\H^g$ obtained by reducing mod $\Gamma_\A$ the
product of the $g$ semi-circle geodesics in $\H$ joining $z_j^-$
to $z_j^+$, $j=1,\ldots,g$. These geodesics are given by the
equations
$$
2\alpha^{(j)}|w_j|^2+\beta^{(j)}(w_j+{\overline
w}_j)+2\gamma^{(j)}=0,\qquad j=1,\ldots,g.
$$
Let ${\mathcal G}_\Delta$ denote the set of representatives of
such $C_h$ mod $\Gamma_\A$. In the case $g=1$, $\Delta=d>0$, these
are the primitive closed geodesics $C\in\Lambda_d$ considered in
Theorem 1 of \cite{Du}.

\goodbreak

\section{Cuspidal Weyl sums and Equidistribution in genus 2}\label{s:5}

As in \S3, let $Q$ be a quadratic form in $m$ variables of
signature $(p,q)$ where $p,q\ge0$ and $p+q=m$ and let $\Gamma$ be
a lattice in $\Omega(Q)$. Let $c_1,c_2,\ldots$ be constants
depending only on $Q$ and $\Gamma$. Let $d\in\Z$, $d\not=0$ be a
square-free integer and $\Gamma_j$, $j=1,\ldots,H(d)$ the
stabilizers in $\Gamma$ of a set of representatives mod $\Gamma$
of the set
$$
\{x\in\R^m: Q[x]=d\}.
$$
Let $\varphi$ be a Maass cusp form of weight 0 on ${\overline{
\Gamma}}\backslash\H_Q$ with eigenvalue $\lambda'$. Then we may
apply Proposition \ref{prop:lift} to
$$
f(z)=\,v^{m/4}<\varphi,{\overline{\theta(z)}}>.
$$
Therefore $f(z)$ is a cuspidal Maass form of discriminant $D$ for
$\Gamma_0(N)\subset\SL(2,\Z)$, where $D$ and $N$ are determined
by $Q$ as in \cite{Du}, p81, and of weight $k=p-\frac m2$ with
eigenvalue $\lambda=\frac14(\lambda'+m-\frac{m^2}4)$. Let
$$
W_{\rm{cusp}}(d,\lambda)=\mu(d)^{-1}
\left(\sum_{i=1}^{H(d)}\int_{\Gamma_i\backslash\Gamma_i(\R)}\varphi(\gamma)d\gamma\right)
$$
where we define
$$
\mu(d)=\sum_{i=1}^{H(d)}{\rm{Vol}}\left(\Gamma_i\backslash\Gamma_i(\R)\right).
$$
From Proposition \ref{prop:weylmaass} we have,
$$
W_{\rm{cusp}}(d,\lambda)=c_3|d|^{m/4}\mu(d)^{-1}\rho(d).
$$
By Siegel's mass formula $\mu(d)$ is a product of local densities.
For $m\ge4$ we have the effective lower bound
$$
\mu(d)\ge c_4|d|^{\frac m2-1}.
$$
Therefore,
$$
W_{\rm{cusp}}(d,\lambda)\le c_5|d|^{1-\frac m4}\rho(d).
$$
By well-known bounds for Fourier coefficients of cusp forms of
integral and half-integral weight (or using the stronger Theorem 5
of \cite{Du}), we know that for $m\ge4$
$$
\lim_{|d|\rightarrow\infty} |d|^{1-\frac m4}\rho(d)=0.
$$
This implies the following result.

\goodbreak

\begin{proposition}\label{weylcusplim}
For $m\ge4$,
$$
\lim_{|d|\rightarrow\infty} W_{\rm{cusp}}(d,\lambda)=0.
$$
\end{proposition}

The Equidistribution in genus 2 statement of \S\ref{s:1} is a
direct corollary of Proposition \ref{weylcusplim} and the
discussion of \S\ref{s:4} once we treat in a similar way the
eigenfunctions of the continuous spectrum of $\Delta_Q$ and show
vanishing results for their Weyl sums. We hope to return to this
in a later paper.

\goodbreak

\section{Cuspidal Weyl sums in the Hilbert modular
case}\label{s:6}

In \S\ref{s:3} the integral quadratic form $Q$ is in $m=3$
variables and is of signature $(p,q)=(2,1)$, and we work over a
totally real field $F$ of degree $g$ over $\Q$ and with the
lattice $\Gamma_\A$ in $G=\SL(2,\R)^g$. Recall that $\A$ is a
fractional ideal in $F$ and that we denote by $\O$ the ring of
integers of $F$ and by $\L$ the lattice $\A^{-1}\oplus\O\oplus\A$
in $F^3$. Let $\Delta\in\O,$ $\Delta\not=0$ and let $\Gamma_i$,
$i=1,\ldots,h(\Delta)$, be the stabilizers in $\Gamma_\A$ of a set
of representatives mod $\Gamma_\A$ of the set
$$
\{h\in\L: Q[h]=\Delta\}.
$$
Let $\varphi\in U$ be an eigenfunction of $\Delta_0^{(j)}$ with
eigenvalue $\lambda_j$, $j=1,\ldots,g$. Then we may apply
Proposition \ref{prop:lifthilbert} to
$$
f(z)=\int_{\Gamma_\A\backslash G}\varphi(g)\theta(z,g)dg.
$$
Suppose $\Delta\ll0$. With the notations of \S\ref{s:3}, Case (i),
let
\begin{equation}\label{weylcuspn}
W_{\rm{cusp}}(\Delta,\lambda)=h(\Delta)^{-1}\sum_{i=1}^{h(\Delta)}\frac1{|\Gamma_i|}\varphi({\overline
g}^{(i)}).
\end{equation}
From the discussion of \S5, we have
\begin{equation}\label{weylcuspnb}
W_{\rm{cusp}}(\Delta,\lambda)=h(\Delta)^{-1}\sum_{z\in\Lambda_{\Delta}}\frac1{|\Gamma_z|}\varphi(z)
\end{equation}
where $\Gamma_z$ is the stabilizer of $z$ in $\Gamma_\A$ and
$h(\Delta)={\rm{Card}}(\Lambda_{\Delta})$. We may take
$$
{\overline
g}^{(i)}_j=\begin{pmatrix}1&h_2^{(i)}/2h_1^{(i)}\\0&1\end{pmatrix}
\begin{pmatrix}\left(\frac{\sqrt{|\Delta_j|}}{2h_1^{(i)}}\right)^{1/2}&0\\0&
\left(\frac{\sqrt{|\Delta_j|}}{2h_1^{(i)}}\right)^{-1/2}\end{pmatrix}.
$$
Then, from Proposition \ref{prop:rhodaven}, we have the following,

\goodbreak

\begin{proposition}\label{prop:weylcuspn} For $\Delta\in\O$, $\Delta\ll0$ we have
$$
W_{\rm{cusp}}(\Delta,\lambda)=2^g(4\pi)^{3g/4}|N_{F/\Q}(\Delta)|^{3/4}h(\Delta)^{-1}\rho(\Delta).
$$
\end{proposition}

Suppose $\Delta\gg0$. With notations as in \S\ref{s:3} let
\begin{equation}\label{weylcuspp}
W_{\rm{cusp}}(\Delta,\lambda)=
\mu(\Delta)^{-1}\sum_{i=1}^{h(\Delta)}\int_{\Gamma_i\backslash\Gamma_i(\R)}\varphi(\gamma)d\gamma,
\end{equation}
where
$$
\mu(\Delta)=\sum_{i=1}^{h(\Delta)}\int_{\Gamma_i\backslash\Gamma_i(\R)}d\gamma.
$$
From the discussion of \S5, we have
\begin{equation}\label{weylcusppb}
W_{\rm{cusp}}(\Delta,\lambda)=\mu(\Delta)^{-1}\sum_{C\in{\mathcal
G}_\Delta}\int_C\varphi(z)ds,
\end{equation}
where
$$
\mu(\Delta)=\sum_{C\in{\mathcal G}_\Delta}{\rm{Vol}}(C).
$$
To explain this last quantity, we consider the Euclidean
hyperbolic distance in each copy of $\H$ in $\H^g$ given by
$$
ds_j^2=y_j^{-2}\left((dx_j)^2+(dy_j)^2\right),\qquad j=1,\ldots,g.
$$
Then we have the real $g$-form
$$
d\underline{s}=\prod_{j=1}^gds_j
$$
and
$$
{\rm{Vol}}(C)=\int_Cd\underline{s}.
$$
By (\ref{rhodavep}) we have the following,

\goodbreak

\begin{proposition}\label{prop:weylcuspp} For $\Delta\in\O$,
$\Delta\gg0$ we have
$$
W_{\rm{cusp}}(\Delta,\lambda)=2^g(4\pi)^{g/4}N_{F/\Q}(\Delta)^{3/4}\mu(\Delta)^{-1}\rho(\Delta).
$$
\end{proposition}

These results enable us bound the cuspidal Weyl sums from above in
terms of the Fourier coefficients of Maass cusp eigenforms of
weight $1/2$ and level $4$. From Propositions \ref{prop:weylcuspn}
and Proposition \ref{prop:weylcuspp} we deduce directly the
following.

\goodbreak

\begin{lemma}\label{lemma:cuspn}
For $\Delta\ll0$ and as $|N_{F/\Q}(\Delta)|\rightarrow\infty$ we
have,
\begin{equation}\label{cuspinn}
|W_{{\mathrm{cusp}}}(\Delta, \lambda)|\ll
\frac{|N_{F/\Q}(\Delta)|^{\frac34}}{h(\Delta)}\rho(\Delta).
\end{equation}
For $\Delta\gg0$ and as $N_{F/\Q}(\Delta)\rightarrow\infty$ we
have,
\begin{equation}\label{cuspinp}
|W_{{\mathrm{cusp}}}(\Delta, \lambda)|\ll
\frac{N_{F/\Q}(\Delta)^{\frac34}}{\mu(\Delta)}\rho(\Delta).
\end{equation}
The implied constants depend only on the field $F$ (and in fact
can be bounded above explicitly by a function of $g$ only).
\end{lemma}

\goodbreak

\section{Eisenstein Weyl sums: the case of the Hilbert modular group for class number 1}\label{s:7}

We continue with the notations of \S\S\ref{s:3}--\ref{s:5}. The
methods of Maass, in particular Proposition \ref{prop:lift} and
Proposition \ref{prop:lifthilbert}, do not apply to the
eigenfunctions of the continuous spectrum of the Laplacian, which
is non-trivial in all cases considered in this paper as the group
actions are not co-compact. These eigenfunctions are furnished by
the Eisenstein series. We are therefore led to consider averages
or Weyl sums as in \S\ref{s:5} with cusp forms replaced by
Eisenstein series.

We shall restrict ourselves to the case $m=3$, $(p,q)=(2,1)$ and
to the Hilbert modular group
$\Gamma=\Gamma_{\A}={\text{PSL}}(2,{\O})$ where ${\O}={\O}_{F}$ is
the ring of integers of a totally real field $F$ of degree $g\ge1$
with class number 1. We repeatedly use facts about Eisenstein
series from \cite{Efr1}, \cite{Shi3}, \cite{Shi4}, \cite{Shi4a},
and \cite{Shi4b}. The generalization to arbitrary class number and
to arbitrary $\Gamma_{\A}$ should be straightforward if somewhat
technical, with some necessary material available in \cite{Sor}.
In order to extend to the case $g>1$ the arguments of \cite{Du} on
Eisenstein Weyl sums, we will generalize in this section some
classical arguments due to Hecke \cite{He} and Kronecker
\cite{Kro}. It would be of interest to consider also the cases
$m=4$, $(p,q)=(2,2)$ and $m=5$, $(p,q)=(3,2)$, in order, for
example, to deduce equidistribution results on the non-cuspidal
part of the corresponding $L^2$-spaces.

The Eisenstein series are $\Gamma$-automorphic eigenfunctions of
the Laplacians $\Delta_0^{(j)}$, $j=1,\ldots,g$ of (\ref{lapz})
corresponding to the continuous part of the spectrum. As we shall
see shortly, they are functions of $(z,s)\in\H^g\times\C$, and
$m\in\Z^{g-1}$. Let,
\begin{multline}\label{lt}
{\rm L}^2(\Gamma\backslash \H^g) =\{\varphi:\H^g\rightarrow\C:
\varphi(\gamma z)=\varphi(z),\,\gamma\in\Gamma,\\
\int_{\Gamma\backslash
\H^g}|\varphi|^2N(v)^{-2}N(du)N(dv)<\infty\}.
\end{multline}
Using general results and more specifically those of \cite{Efr1}
and \cite{Shi4a}, we have a decomposition,
$$
{\rm L}^2(\Gamma\backslash \H^g)={\rm
L}^2_{\rm{cusp}}(\Gamma\backslash \H^g)\oplus {\mathcal R}\oplus
{\mathcal E},$$ where ${\mathcal E}$ is generated in an
appropriate $L^2$-sense by the Eisenstein series evaluated at
$s=\frac12+it$, $t\in\R$, and ${\mathcal R}$ is generated by the
residues of their finitely many poles in $s\in(\frac12, 1]$ (the
Eisenstein series are not themselves $L^2$-integrable). For
$g\ge1$, the only such pole occurs at $s=1$ and has residue given
by the volume of the fundamental region of $\Gamma$. The cuspidal
part of the $L^2$-decomposition is given by (\ref{ltu}).

As $F$ is assumed to have class number 1, the group $\Gamma$ has 1
cusp at infinity with stabilizer $\Gamma_{\infty}$. The Eisenstein
series are of the form
\begin{equation}\label{eisdef}
E(z,s,m)=\sum_{\gamma\in\Gamma_{\infty}\backslash\Gamma}y^s(\gamma(z))\lambda_m(y(\gamma(z))),\qquad
{\rm{Re}}(s)>1.
\end{equation}
Here $z=(z_j)_{j=1}^g\in\C^g$ with $z_j=x_j+iy_j$, $x_j,y_j\in\R$
and $y(\gamma(z))=\prod_{j=1}^g{\rm{Im}}(\gamma^{(j)}(z_j))$ where
$\gamma^{(j)}$ is obtained from $\gamma$ by applying the $j$-th
Galois embedding to its entries, for a given ordering of the
Galois embeddings of $F$ into $\R$, starting with the identity
embedding. Moreover, we have $s\in\C$ and
$m=(m_q)_{q=1}^{g-1}\in\Z^{g-1}$ with $\lambda_m$ an exponential
sum similar to a Gr\"ossencharacter. Namely, there are parameters
$e_j^{(q)}\in\R$, $q=1,\ldots,g-1$, $j=1,\ldots,g$, determined by
the choice of basis of the unit group $\O^\ast$ of $\O$, such that
$$
\lambda_m(z)=\prod_{j=1}^g\prod_{q=1}^{g-1}|z_j|^{\pi i
m_qe_j^{(q)}}.
$$
Moreover, if
$\varepsilon_q=(\varepsilon_q^{(1)},\ldots,\varepsilon_q^{(g)})$,
$q=1,\ldots,g-1$, are totally positive generators of the group
$\O^\ast$ embedded in ${(\R^\ast)}^g$ by the Galois embeddings of
$F$ into $\R$, then the $e_j^{(q)}$ are determined by the
equations
$$
\sum_{j=1}^ge_j^{(q)}=0,\qquad q=1,\ldots,g-1,$$ and
$$\sum_{j=1}^ge_j^{(q)}\log\varepsilon_r^{(j)}=\delta_{r,q},\qquad
r,q=1,\ldots,g-1.
$$
These conditions ensure, in particular, that for
$u=(u^{(1)},\ldots,u^{(g)})$ in $\O^\ast$ we have,
\begin{equation}\label{unit}
\prod_{j=1}^g\prod_{q=1}^{g-1}\mid{u^{(j)}}\mid^{\pi i
m_qe_j^{(q)}}=1.
\end{equation}
The series $E(z,s,m)$ has a meromorphic continuation to all of
$s\in\C$. Full details can be found in \cite{Efr1}, Chapter II. We
shall often perform formal manipulations with the series
definition of $E(z,s,m)$ without specifying each time the domain
of convergence.

We have, using \cite{Efr1}, p47,
\begin{equation}\label{eisredef}
E(z,s,m)=\sum_{\{c,d\},
(c,d)=1}\prod_{j=1}^g\frac{y_j^{s_j}}{|c^{(j)}z_j+d^{(j)}|^{2s_j}}
\end{equation}
where the summation is over $c,d\in\O$, with $\{c,d\}$ meaning
that pairs differing by multiplication by an element of
$\O^{\ast}$ are identified and $(c,d)=1$ means that $c,d$ generate
$\O$. As usual $c^{(j)}$, $d^{(j)}$ are the $j$-th Galois
conjugates of $c$, $d$. Furthermore,
\begin{equation}s_j=s+\pi i\sum_{q=1}^{g-1}m_qe_j^{(q)}, \qquad j=1,\ldots,g.
\end{equation}
In that same reference, it is shown that if
\begin{equation}\label{eisf}
F(z,s,m)=\sum_{\{c,d\}}\prod_{j=1}^g\frac{y_j^{s_j}}{|c^{(j)}z_j+d^{(j)}|^{2s_j}}
\end{equation}
then,
\begin{equation}\label{eisfid}
F(z,s,m)=L(2s,\lambda_{-2m})E(z,s,m),
\end{equation}
where $L(s,\lambda_m)$ is the Hecke zeta function given by
\begin{equation}\label{lgr}
L(s,\lambda_m)=\sum_{(b)}\frac{\lambda_m(b)}{|N_{F/\Q}(b)|^s}.
\end{equation}
Here, the sum is over the (principal) integral ideals $(b)$ of $F$
and
$$
\lambda_m(b)=\prod_{j=1}^g\prod_{q=1}^{g-1}|b^{(j)}|^{\pi i
m_qe_j^{(q)}}.
$$
This last expression is well-defined thanks to (\ref{unit}).

With the notation of \S\ref{s:4} and \S\ref{s:5}, the Eisenstein
Weyl sums are given, for $\Delta\ll0$, by
\begin{equation}\label{weyleisn}
W_{\rm{Eis}}(\Delta,t,m)=\frac1{h(\Delta)}\sum_{z_h\in\Lambda_{\Delta}}E(z_h,\frac12+it,m),
\end{equation}
where $h(\Delta)$ is the cardinality of $\Lambda_{\Delta}$ and,
for $\Delta\gg0$, by
\begin{equation}\label{weyleisp}
W_{\rm{Eis}}(\Delta,t,m)=\frac1{\mu(\Delta)}\sum_{C\in{\mathcal
G}_\Delta}\int_CE(z,\frac12+it,m)d\underline{s},
\end{equation}
where
$$
\mu(\Delta)=\sum_{C\in{\mathcal G}_\Delta}{\rm{Vol}}(C).
$$
In order to relate these Weyl sums to Fourier coefficients of
half-integral weight Hilbert-Maass Eisenstein series, we need to
generalize for $g>1$ some classical arguments due to Hecke
\cite{He} and Kronecker \cite{Kro} that apply to the special case
$g=1$ as in \cite{Du}. In particular, the formulae of Proposition
\ref{prop:weyln} and Proposition \ref{prop:weylp} at the end of
this section generalize these classical results.

From now on, we assume that $\Delta\in\O$ generates the ideal
given by the relative discriminant of $L=F(\sqrt{\Delta})$ over
$F$. This replaces the fundamental discriminant assumption in
\cite{Du}, Theorem 1. As $F$ has class number 1, the ideals of $L$
are free $\O$-modules of rank 2. Let
$\rho\rightarrow{\overline{\rho}}$, $\rho\in L$, denote the
non-trivial automorphism of $L$ over $F$. We may choose a relative
basis $\{1,\Omega\}$ of $\O_L$ over $\O$ such that
$\Delta=\left({\overline\Omega}-\Omega\right)^2.$ Let
${\O}^\ast_L$ denote the units of $\O_L$ and $\chi_{L/F}$ denote
the relative field character for $L$ over $F$.

Consider first the case $\Delta\ll0$. By (\ref{eisfid}) we have
\begin{equation}\label{weyln2}
L(2s,\lambda_{-2m})\sum_{z_h\in\Lambda_{\Delta}}E(z_h,s,m)=
\sum_{z_h\in\Lambda_{\Delta}}\sum_{\{c,d\}}\prod_{j=1}^gQ_h^{(j)}(c^{(j)},d^{(j)})^{-s_j},
\end{equation}
where
\begin{equation}\label{weyln3}
Q_h^{(j)}(c^{(j)},d^{(j)})=\frac{|c^{(j)}z_j+d^{(j)}|^2}{y_j},\qquad
j=1,\ldots,g.
\end{equation}
In the notation of \S\ref{s:4} (except we denote $z_j^+$ by
$z_j$), for $h=(\alpha,\beta,\gamma)\in\O^3$, we have
\begin{equation}\label{weyln4}
Q_h^{(j)}(c^{(j)},d^{(j)})=\frac{2}{\sqrt{|\Delta^{(j)}|}}\;q_h^{(j)}(-d^{(j)},c^{(j)}),
\end{equation}
where
\begin{equation}\label{weyln5}
q_h^{(j)}(x,y)=\alpha^{(j)}x^2+\beta^{(j)}xy+\gamma^{(j)}y^2=\alpha^{(j)}(x-z_jy)(x-{\overline
z_j}y).
\end{equation}
Therefore,
\begin{equation}\label{weyln6}
\prod_{j=1}^gQ_h^{(j)}(c^{(j)},d^{(j)})^{-s_j}=\prod_{j=1}^g|\Delta^{(j)}/4|^{s_j/2}\prod_{j=1}^g
q_{h}^{(j)}(-d^{(j)},c^{(j)})^{-s_j}.
\end{equation}
There is a bijection between the ideal classes of $\O_L$ and the
points $z_h\in\Lambda_{\Delta}$. To $h=(\alpha,\beta,\gamma)\in
W_{\Delta}$, we associate the ideal ${\mathcal A}_h$ with basis
$\{\alpha,\frac{\beta+\sqrt\Delta}2\}$. The relative norm over $F$
of ${\mathcal A}_h$ is generated by $\alpha$. The relative norms
of the integral ideals in the ideal class ${\rm{Cl}}({\mathcal
A}_h)$ of ${\mathcal A}_h$ are generated by $q_h(-d,c)$, for
$c,d\in\O$. Therefore,
\begin{multline}\label{weyln7}
\sum_{\{c,d\}}\prod_{j=1}^gQ_h^{(j)}(c^{(j)},d^{(j)})^{-s_j}=2^{-sg}|N_{F/\Q}(\Delta)|^{s/2}\lambda_{m/2}(\Delta/4)\times\cr\times\sum_{{\mathcal
A}\in{\rm{Cl}({\mathcal
A}_h)}}\frac{\lambda_{-m}(N_{L/F}({\mathcal
A}))}{N_{L/\Q}({\mathcal A})^s}.
\end{multline}
Combining (\ref{weyln2}) and (\ref{weyln7}) we conclude that
\begin{multline}\label{weyln8}
2^{sg}L(2s,\lambda_{-2m})\sum_{z_h\in\Lambda_{\Delta}}E(z_h,s,m)=|N_{F/\Q}(\Delta)|^{s/2}\lambda_{m/2}(\Delta/4)\times\cr\times
L(s,\lambda_{-m},L),
\end{multline}
where
\begin{equation}\label{weyln9}
L(s,\lambda_{-m},L)=\sum_{\mathcal
A}\frac{\lambda_{-m}(N_{L/F}({\mathcal A}))}{N_{L/\Q}({\mathcal
A})^s}=\sum_{\mathsf A}\sum_{{\mathcal A}\in{\mathsf
A}}\frac{\lambda_{-m}(N_{L/F}({\mathcal A}))}{N_{L/\Q}({\mathcal
A})^s},
\end{equation}
with ${\mathcal A}$ ranging over the non-zero integral ideals of
$L$ and ${\mathsf A}$ ranging over the ideal classes of $L$. From
(\ref{weyln8}) and (\ref{weyln9}) we deduce,
\begin{multline}\label{weyln10}
2^{sg}L(2s,\lambda_{-2m})\sum_{z_h\in\Lambda_{\Delta}}E(z_h,s,m)=|N_{F/\Q}(\Delta)|^{s/2}\lambda_{m/2}(\Delta/4)\times\cr\times
L(s,\lambda_{-m})L(s,\chi_{L/F}\lambda_{-m}),
\end{multline}
where
\begin{equation}\label{weyln11}
L(s,\chi_{L/F}\lambda_{-m})=\sum_{(b)}\frac{\chi_{L/F}(b)\lambda_{-m}(b)}{|N_{F/\Q}(b)|^s},
\end{equation}
the sum ranging over the ideals $(b)$ of $F$. We deduce finally
the following result.

\goodbreak

\begin{proposition}\label{prop:weyln} For $\Delta\ll0$ we have,
\begin{multline}\label{weyln12}
W_{\rm{Eis}}(\Delta,t,m)=2^{-sg}\frac{L(\frac12+it,\lambda_{-m})}{L(1+it,\lambda_{-2m})}\times\cr\times
\frac{|N_{F/\Q}(\Delta)|^{\frac14+it}\lambda_{m/2}(\Delta/4)}{h(\Delta)}L(\frac12+it,\chi_{F(\sqrt{\Delta})/F}\lambda_{-m}).
\end{multline}
\end{proposition}

Now consider the case $\Delta\gg0$. Combining (\ref{eisf}) and
(\ref{weyleisp}) we have
\begin{equation}\label{weylp1}
L(2s,\lambda_{-2m})\sum_{C_h\in{\mathcal
G}_\Delta}\int_{C_h}E(z,s,m)d{\underline{s}}=
\sum_{C_h\in{\mathcal
G}_\Delta}\int_{C_h}\sum_{\{c,d\}}\prod_{j=1}^g\frac{y_j^{s_j}}{|c^{(j)}z+d^{(j)}|^{2s_j}}d\underline{s}.
\end{equation}
We adapt some ideas of Hecke \cite{He} for the case $g=1$, as they
are explained in \cite{Za1}, \cite{Za2}, to the general case
$g\ge1$. Let ${\mathsf A}$ be a fixed ideal class of
$L=F(\sqrt\Delta)$ and ${\mathcal B}$ a fixed element of ${\mathsf
A}^{-1}$. We have a correspondence ${\mathcal A}\mapsto {\mathcal
A}{\mathcal B}=(\eta)$, which is a bijection between the set of
ideals of ${\mathsf A}$ and the set of principal ideals with
$\eta\in{\mathcal B}$.

For a fixed ideal class ${\mathsf A}$, define
\begin{equation}\label{weylp2}
L(s,\lambda_{-m},{\mathsf A})=\sum_{{\mathcal A}\in{\mathsf
A}}\frac{\lambda_{-m}(N_{L/F}({\mathcal A}))}{N_{L/Q}({\mathcal
A})^s}.
\end{equation}
Define as in (\ref{weyln9}),
\begin{equation}\label{weylp3}
L(s,\lambda_{-m},L)=\sum_{\mathcal
A}\frac{\lambda_{-m}(N_{L/F}({\mathcal A}))}{N_{L/Q}({\mathcal
A})^s}=\sum_{\mathsf A}L(s,\lambda_{-m},{\mathsf A}),
\end{equation}
with ${\mathcal A}$ ranging over the non-zero integral ideals of
$L$ and ${\mathsf A}$ ranging over the ideal classes of $L$. For
${\mathcal B}$ a fixed element of ${\mathsf A}^{-1}$, we have
\begin{equation}\label{weylp4}
L(s,\lambda_{-m},{\mathsf A})=N_{L/\Q}({\mathcal B})^s
\lambda_m(N_{L/F}({\mathcal B}))\sum_{{\mathcal A}\in{\mathsf
A}}\frac{\lambda_{-m}(N_{L/F}({\mathcal AB}))}{N_{L/\Q}({\mathcal
AB})^s}.
\end{equation}
Two numbers $\eta_1$, $\eta_2\in{\mathcal B}$ define the same
principal ideal if and only if $\eta_1=\varepsilon\eta_2$, for
$\varepsilon\in{\O}^\ast_L$. Hence,
\begin{equation}\label{weylp5}
L(s,\lambda_{-m},{\mathsf A})=N_{L/\Q}({\mathcal B})^s
\lambda_m(N_{L/F}({\mathcal B})){\sum}'_{\eta\in{\mathcal
B}/\O^\ast_L}\frac{\lambda_{-m}(N_{L/F}(\eta))}{|N_{L/\Q}(\eta)|^s},
\end{equation}
with $\sum'$ denoting a sum over non-zero elements. We have an
exact sequence
\begin{equation}\label{weylp6}
1\rightarrow\O_{L,1}^\ast\rightarrow\O_L^\ast\rightarrow\O^\ast,
\end{equation}
where the right most arrow is given by the reduced norm from $L$
to $F$ and $\O_{L,1}^\ast$ is the group of units of $\O_L$ of
reduced norm 1. The image $N_{L/F}(\O_L^\ast)$ is of finite index
in $\O^\ast$. Moreover, $\O_{L,1}^\ast$ is a free abelian group of
rank $g$ (as remarked already in \S\ref{s:3}). Notice that
$\O_{L,1}^\ast\cap\O^\ast=\{\pm1\}$ and that the group
$\O_{L,1}^\ast\O^\ast$ is of finite index $\mathsf i$ in
$\O^\ast_L$. We have therefore from (\ref{weylp5}),
\begin{equation}\label{weylp7}
L(s,\lambda_{-m},{\mathsf A})=N_{L/\Q}({\mathcal B})^s
\lambda_m(N_{L/F}({\mathcal B}){\mathsf
i}^{-1}S(s,\lambda_{-m},{\mathcal B}),
\end{equation}
where
\begin{equation}\label{weylp8}
S(s,\lambda_{-m},{\mathcal B})={\sum}'_{\eta\in{\mathcal
B}/\O_{L,1}^\ast\O^\ast}\frac{\lambda_{-m}(\eta{\overline
{\eta}})}{|N_{L/\Q}(\eta)|^s}.
\end{equation}
Let $\varepsilon^{(i)}$, $i=1,\ldots,g$, be generators of
$\O^\ast_{L,1}/\{\pm1\}$. Let $\xi_j$ be the extension to $L$ of
the $j$-th Galois embedding of $F$ into $\R$ chosen so that
$\xi_j(\sqrt\Delta)=\sqrt{\Delta^{(j)}}>0$, $j=1,\ldots,g$. Let
$\eta_j=\xi_j(\eta)$ for $\eta\in L$. We may suppose that
$\varepsilon^{(i)}_j>0$ for $i,j=1,\ldots,g$. We have
\begin{equation}\label{weylp9}
S(s,\lambda_{-m},{\mathcal B})={\sum}'_{\eta\in{\mathcal
B}/\O^\ast\langle\varepsilon^{(1)},\ldots,\varepsilon^{(g)}\rangle}|\eta_1{\overline{\eta_1}}|^{-s_1}\ldots|\eta_g{\overline{\eta_g}}|^{-s_g}.
\end{equation}
Hecke observed the following identity (see \cite{Za1}, p161): for
$a,b\in\R$, $a,b\not=0$,
\begin{equation}\label{weylp10}
\int_{-\infty}^\infty\frac{dv}{(a^2e^v+b^2e^{-v})^s}=\frac{c(s)}{|ab|^s},
\end{equation}
where
\begin{equation}\label{weylp11}
c(s)=\int_{-\infty}^\infty\frac{dv}{(e^v+e^{-v})^s}.
\end{equation}
Let $N(c(s))=\prod_{j=1}^gc(s_j)$. We deduce that,
\begin{equation}\label{weylp12}
N(c(s))S(s,\lambda_{-m},{\mathcal B})={\sum}'_{\eta\in{\mathcal
B}/\O^\ast\langle\varepsilon^{(1)},\ldots,\varepsilon^{(g)}\rangle}
\prod_{j=1}^g\int_{-\infty}^\infty\frac{dv_j}{(\eta_j^2e^{v_j}+{\overline{\eta_j}}^2e^{-v_j})^{s_j}}.
\end{equation}
Using the embeddings $\xi_j$,$j=1,\ldots,g$, we may embed any
element of $L$ into $\R^g$. The transformation
$\eta\mapsto\varepsilon\eta$, $\varepsilon\in\O_{L,1}^\ast$ with
$\varepsilon_j>0$, $j=1,\ldots,g$, then corresponds to a vector
translation given componentwise by
$$v_j\mapsto v_j+2\log\varepsilon_j,\qquad j=1,\ldots,g.$$
Let ${\mathcal L}_g$ be the lattice in $\R^g$ generated by the
vectors $(2\log\varepsilon^{(i)}_j)_{j=1}^g$, $i=1,\ldots,g$.
Then, from (\ref{weylp12}) we deduce that,
\begin{equation}\label{weylp13}
N(c(s))S(s,\lambda_{-m},{\mathcal B})={\sum}'_{\eta\in{\mathcal
B}/\O^\ast}\int_{{\mathcal
L}_g\backslash\R^g}\frac{N(dv)}{\prod_{j=1}^g(\eta_j^2e^{v_j}+{\overline{\eta}}_j^2e^{-v_j})^{s_j}},
\end{equation}
where $N(dv)=\prod_{j=1}^gdv_j$.

Now suppose the ideal ${\mathcal B}$ has basis $\{1,w\}$ (with
$w>{\overline{w}}$). Then $\eta=cw+d$ for $c,d\in\O$ and
$\eta_j=c^{(j)}w_j+d^{(j)}$ and,
\begin{equation}\label{weylp14}
\eta_j^2e^{v_j}+{\overline{\eta}}_j^2e^{-v_j}=(c^{(j)}w_j+d^{(j)})^2e^{v_j}+(c^{(j)}{\overline
w}_j+d^{(j)})^2e^{-v_j}.
\end{equation}
Let $w^+_j=\max(w_j,{\overline{w}}_j)$ and
$w^-_j=\min(w_j,{\overline{w}}_j), j=1,\ldots,g$. Make the change
of variables,
\begin{equation}\label{weylp15}
z_j=\frac{w_j^+\sqrt{-1}e^{v_j}+w^-_j}{\sqrt{-1}e^{v_j}+1},\qquad
j=1,\ldots,g.
\end{equation}
Then as $v_j$ ranges from $-\infty$ to $\infty$ the variable $z_j$
runs over the geodesic in $\H$ joining $w^-_j$ to $w^+_j$. A
direct calculation shows, with $y_j={\mathrm{Im}}(z_j)$, that
\begin{equation}\label{weylp16}
y_j|c^{(j)}z_j+d^{(j)}|^{-2}=(w^+_j-w^-_j)\lbrace
e^{v_j}(c^{(j)}w_j+d^{(j)})^2+e^{-v_j}(c^{(j)}w_j+d^{(j)})^2\rbrace^{-1}.
\end{equation}
In the notation of \S\ref{s:4}, let $h=(\alpha,\beta,\gamma)\in
W_{\Delta}$ and let $w_j^+=z_j^+$, $w^-_j=z_j^{-}$. Then
${\mathcal B}={\mathcal B}_h=\O+\O z_1^+$ and
$w_j^+-w^-_j=\sqrt{\Delta^{(j)}}/|\alpha^{(j)}|$ if
$\alpha\not=0$. Under the change of variables (\ref{weylp15}), the
quotient ${\mathcal L}_g\backslash\R^g$ becomes the quotient
$\Gamma_h\backslash\Gamma_h(\R)$ realized as ${\mathcal C}_h$.
From (\ref{weylp12}) and (\ref{weylp16}) we deduce
\begin{multline}\label{weylp17}
N(c(s))S(s,\lambda_{-m},{\mathcal
B})=\prod_{j=1}^g\left(\sqrt{\Delta^{(j)}}/|\alpha^{(j)}|\right)^{-s_j}\times\cr\times\sum_{\{c,d\}}\int_{{\mathcal
C}_h}\prod_{j=1}^gy_j^{s_j}|c^{(j)}z_j+d^{(j)}|^{-2s_j}d{\underline{s}}.
\end{multline}
Now,
\begin{equation}\label{weylp18}
N_{L/\Q}({\mathcal B})^s \lambda_m(N_{L/F}({\mathcal
B}))=\prod_{j=1}^g\left(N_{L/F}({\mathcal
B})\right)^{s_j}=\prod_{j=1}^g|\alpha^{(j)}|^{-s_j}.
\end{equation}
Let ${\mathsf A}_h$ denote the ideal class of ${\mathcal B}_h$.
Combining (\ref{eisf}), (\ref{eisfid}), (\ref{weylp7}),
(\ref{weylp8}), (\ref{weylp17}), and (\ref{weylp18}) we deduce
that,
\begin{equation}\label{weylp19}
L(2s,\lambda_{-2m})\int_{{\mathcal
C}_h}E(z,s,m)d{\underline{s}}=N(c(s)){\mathsf
i}N_{F/\Q}(\Delta)^{s/2}\lambda_{m/2}(\Delta)L(s,\lambda_{-m},{\mathsf
A}_h).
\end{equation}
There is a bijection between the ideal classes of $\O_L$ and the
representatives of ${\mathcal G}_{\Delta}$ with ${\mathcal C}_h$
corresponding to ${\mathsf A}_h$, $h\in W_{\Delta}$. We conclude
that
\begin{equation}\label{weylp20}
L(2s,\lambda_{-2m})\sum_{{\mathcal C}_h\in{\mathcal
G}_{\Delta}}\int_{{\mathcal
C}_h}E(z,s,m)d{\underline{s}}=N(c(s)){\mathsf
i}N_{F/\Q}(\Delta)^{s/2}\lambda_{m/2}(\Delta)L(s,\lambda_{-m},L),
\end{equation}
and finally that
\begin{multline}\label{weylp21}
L(2s,\lambda_{-2m})\sum_{{\mathcal C}_h\in{\mathcal
G}_{\Delta}}\int_{{\mathcal
C}_h}E(z,s,m)d{\underline{s}}=N(c(s)){\mathsf
i}N_{F/\Q}(\Delta)^{s/2}\lambda_{m/2}(\Delta)\times\cr \times
L(s,\lambda_{-m})L(s,\chi_{L/F}\lambda_{-m}).
\end{multline}
We deduce finally the following result.

\goodbreak

\begin{proposition}\label{prop:weylp} For $\Delta\gg0$ we have,
\begin{multline}\label{weylp22}
W_{{\mathrm{Eis}}}(\Delta,t,m)=N(c(s)){\mathsf
i}\frac{L(\frac12+it,\lambda_{-m})}{L(1+it,\lambda_{-2m})}\times\cr
\times
\frac{N_{F/\Q}(\Delta)^{\frac14+it}\lambda_{m/2}(\Delta)}{\mu(\Delta)}L(\frac12+it,\chi_{F(\sqrt\Delta)/F}\lambda_{-m}).
\end{multline}
\end{proposition}

From Propositions \ref{prop:weyln} and Proposition
\ref{prop:weylp} we deduce directly the following.

\goodbreak

\begin{lemma}\label{lemma:eisenweyl}
For $\Delta\ll0$ and as $|N_{F/\Q}(\Delta)|\rightarrow\infty$ we
have,
\begin{equation}\label{eisenp}
|W_{{\mathrm{Eis}}}(\Delta, t,m)|\ll
\frac{|N_{F/\Q}(\Delta)|^{\frac14}}{h(\Delta)}L(\frac12+it,\chi_{F(\sqrt\Delta)/F}\lambda_{-m}).
\end{equation}
For $\Delta\gg0$ and as $N_{F/\Q}(\Delta)\rightarrow\infty$ we
have,
\begin{equation}\label{eisenn}
|W_{{\mathrm{Eis}}}(\Delta, t,m)|\ll
\frac{N_{F/\Q}(\Delta)^{\frac14}}{\mu(\Delta)}L(\frac12+it,\chi_{F(\sqrt\Delta)/F}\lambda_{-m}).
\end{equation}
The implied constants depend on $t$ and on the field $F$.
\end{lemma}

In \cite{Du}, the Eisenstein Weyl sums for $g=1$ are shown to be
proportional to the Fourier coefficients of Eisenstein series of
half-integral weight and level 4 using explicit formulae for these
coefficients derived in \cite{GH}. In \cite{Shi4}, general
formulae for Fourier coefficients of Eisenstein series of
half-integral weight and level dividing $4$ for the group
${\mathrm{Sp}}(m,F)$ are obtained, where $F$ is a totally real
algebraic number field. The case $m=1$ gives generalizations of
the formulae of \cite{GH} to the Hilbert modular case. In the
notations of \cite{Shi4}, Theorem 6.1, the product
$L(s,\lambda_{-m})L(s,\chi_{F(\sqrt\Delta)/F}\lambda_{-m})$
occurring (at $s=\frac12+it$) in Proposition \ref{prop:weyln} and
Proposition \ref{prop:weylp} is proportional to the product of the
first Fourier coefficient and the ``$\Delta$''-th Fourier
coefficient $c_{\mathbf f}(\Delta,s)$ of the Eisenstein series
$E'=E'(z,s,1/2,0,\lambda_{-m},4)$ of level $4$ and weight $1/2$.
The $c_{\mathbf f}(\Delta,s)$ of \cite{Shi4} correspond to
$|N_{F/\Q}(\Delta)|^{1/2}\rho(\Delta, E')$ with the conventions
that we adopt in \S4 (\ref{fuv}).

These results enable us to bound the Eisenstein Weyl sums from
above in terms of the Fourier coefficients of Eisenstein series of
weight $1/2$ and level $4$, or alternatively central values of
$L$-functions.
\goodbreak

\section{Expected subconvexity results and proof of Theorem \ref{equigrh}}\label{s:8}

We continue with the assumptions and notations of \S\ref{s:6} and
\S\ref{s:7}. The equidistribution results of Theorem \ref{equigrh}
would follow, without GRH, from an unconditional proof of
\begin{equation}\label{weylclim}
\lim_{|N_{F/\Q}(\Delta)|\rightarrow\infty}W_{\mathrm{cusp}}(\Delta,\lambda)=0,
\end{equation}
and
\begin{equation}\label{Eisclim}
\lim_{|N_{F/\Q}(\Delta)|\rightarrow\infty}W_{\mathrm{Eis}}(\Delta,t,m)=0.
\end{equation}
As $\Delta$ is a fundamental (relative) discriminant which is
totally definite, the number $h(\Delta)$ can be replaced by the
class number of $\O_L$, $L=F(\sqrt\Delta)$ (see for example
\cite{Coh}, Chapter 7, \cite{Mas}). Moreover, for $\Delta\gg0$, we
have by \cite{Efr1}, p36, that $\mu(\Delta)=h(\Delta)R$ where $R$
is a regulator associated to $\O^\ast_{L,1}$ and given by
$$
R=\det\left(2\log\varepsilon_j^{(i)}\right)_{i,j=1}^g.
$$
The results of \S\ref{s:6} and \S\ref{s:7} show that the
(ineffective) lower bounds
$$
h(\Delta)\gg_{\varepsilon}|N_{F/\Q}(\Delta)|^{1/2-\varepsilon},\qquad
{\mathrm{as}}\quad |N_{F/\Q}(\Delta)|\rightarrow\infty,
$$
and
$$
h(\Delta)R\gg_{\varepsilon}N_{F/\Q}(\Delta)^{1/2-\varepsilon},\qquad
{\mathrm{as}}\quad N_{F/\Q}(\Delta)\rightarrow\infty,
$$
provided by the Brauer-Siegel Theorem \cite{Bra}, \cite{Sie0}
together with generalizations to the case $g>1$ of the
subconvexity results for $g=1$ in \cite{Du}, Theorem 5 would imply
(\ref{weylclim}) and (\ref{Eisclim}). The corresponding
subconvexity results for the holomorphic case have been shown in
\cite{CoSaPi}. We would need the Fourier coefficients
$\rho(\Delta,f)$ for $f$ a cusp form with $L^2$-norm 1 or an
Eisenstein series, with eigenvalue $\lambda$ and half-integral
weight $k$, to have an upper bound in the $\Delta$-aspect as good
as $\rho(\Delta, f)\ll_{k,\epsilon}
c(\lambda)|N_{F/\Q}(\Delta)|^{-1/4-\delta+\epsilon}$ for a fixed
$\delta>0$ and a positive explicit constant $c(\lambda)$. From
Lemma \ref{lemma:eisenweyl}, the desired result for Eisenstein
series would follow from a subconvexity result for $L$-functions
of ${\rm GL}(1,F)$. Partial progress towards subconvexity results
in the Maass case have been made by Gergely Harcos \cite{Har}, but
the complete adaptation of the ${\rm GL}(2,F)$ methods of
\cite{CoSaPi} to the Maass case remains elusive. Such results
would follow however from GRH, so our Theorem \ref{equigrh}
remains conditional.

\bigskip

\noindent{\bf{Note added in Proof:}} As remarked to us by Emmanuel
Ullmo, some comments must be added about the dependence on the
parameter $m\in\Z^{g-1}$ in the upper bounds for Eisenstein
series. We will add such a comment shortly.

\vfill\break
\bibliographystyle{amsalpha}

\end{document}